\documentclass[12pt, reqno]{amsart}
\usepackage{preemble}
\usepackage{microtype}
\begin{document}

\title[Cutoff for congestion dynamics]
      {Cutoff for congestion dynamics and related generalized exclusion processes}
\author{Ryokichi Tanaka}
\address{Department of Mathematics, Kyoto University, Kyoto 606-8502 JAPAN}
\email{rtanaka@math.kyoto-u.ac.jp}
\date{\today}

\maketitle

\begin{abstract}
We consider congestion dynamics with $n$ players and $Q$ resources under the constraint that the number of each resource is $\kappa$ and that $n<\kappa Q$ in the regime that
 $n$ and $\kappa$ diverge but $Q$ is fixed with $n=\floor{\rho \kappa Q}$ for a fixed constant $\rho \in (0, 1/2]$.
 We show that the Glauber dynamics and its unlabeled version exhibit cutoff at time $(1/2)n \log n$ and $(1/2)(1-\rho)n\log n$ in total variation respectively.
 The unlabeled version is a special case of natural Markov chains for sampling from log M-concave distributions.
 We also show that a family of Markov chains for uniform sampling on M-convex sets does not necessarily exhibit cutoff.
\end{abstract}

\section{Introduction}\label{Sec:introduction}

This article aims to study the following special case of {\bf congestion dynamics}.
Let $n$, $Q$ and $\kappa$ be positive integers.
We consider an evolution of configurations where each one of $n$ players adapts one of $Q$ resources under the constraint that the number of each resource is $\kappa$ and that $n<\kappa Q$.
At each step, one of players chosen uniformly at random releases the current resource and adopt a new one independently with the probability proportional to the number of leftovers.
This model is a special case of Markov chains for sampling from Gibbs distributions (with log-concave potentials) by the Glauber dynamics \cite{logconcaveIV} and \cite{Kleer},
and interacting particle systems \cite{KipnisLandim}.
We will discuss the background in Section \ref{Sec:background} below.

A main focus will be in the regime where the number of players $n$ diverges but $Q$ is fixed and the ratio $n/(\kappa Q)$ is a fixed constant $\rho \in (0, 1/2]$ up to rounding.
We show that this particular dynamics and its unlabeled version (i.e., players are indistinguishable) reveal abrupt convergences to their unique equilibrium states.
More precisely,
they exhibit {\bf cutoff} at time $(1/2)n\log n$ and $(1/2)(1-\rho)n\log n$ as $n \to \infty$
in the worst case total variation distance respectively.
On the one hand, 
the results are expected to hold in a more general setting in the same category we are considering.
On the other hand, 
it turns out that the unlabeled process falls into a class of {\bf generalized exclusion processes} on complete graphs (cf.\ \cite[Section 4, Chapter 2]{KipnisLandim}).
The setting extends exclusion processes on these graphs (or, the Bernoulli-Laplace urn model)
for which a cutoff has been studied in \cite{DiaconisShahshahaniBernoulliLaplace},
\cite{LacoinLeblond} and \cite{PatkoPete}; see also a more recent result \cite{Olesker-Taylor-Schmid}.

Let us state the precise statements.
The set of {\bf resources} is defined by $S:=\{1, \dots, Q\}$.
For $\sigma=(\sigma_i)_{i\in [n]}\in S^n$ where $[n]:=\{1, \dots, n\}$,
let
\[
\xi_v(\sigma):=\#\big\{i \in [n] \ : \ \sigma_i=v\big\} \quad \text{for each $v \in S$}.
\]
In the above $\#$ of a set stands for the cardinality.
The space of {\bf configurations} is 
\[
\Sc_n:=\big\{\sigma \in S^n \ : \ \xi_v(\sigma) \le \kappa \quad \text{for all $v \in S$}\big\}.
\]
This is the state space where an evolution takes place,
and each $\sigma$ indicates an assignment to each player one of $Q$ resources keeping each resource being adapted by at most $\kappa$ players.

Let us fix a real $\rho \in (0, 1)$.
We define $N:=\kappa Q$ and $n=\floor{\rho N}$,
where $\floor{r}$ of a real $r$ is the largest integer at most $r$. 
Let us consider the function $\f:\{0, 1, \dots, \kappa\}\to \R_+$ defined by
\[
\f(s)=\log(\kappa+1)-\log(\kappa+1-s).
\]
The {\bf potential} $\F:\Sc_n\to \R_+$ is defined by $\F(\sigma):=\sum_{v \in S}\sum_{s=0}^{\xi_v(\sigma)}\f(s)$.
The associated {\bf Gibbs distribution} on $\Sc_n$ has the form
\[
\mu_n(\sigma):=\frac{1}{Z_n}e^{-\F(\sigma)} \quad \text{where $Z_n:=\sum_{\sigma \in \Sc_n}e^{-\F(\sigma)}$}.
\]

The process on the state space $\Sc_n$ we are considering evolves by the {\bf Glauber dynamics}:
Given a configuration $\sigma \in \Sc_n$,
first choose $i \in [n]$ uniformly at random,
then replace $\sigma$ by $\sigma'$ according to the Gibbs distribution conditioned that $\sigma'$ coincides with $\sigma$ except for $i$.
In this special form of function $\f$, 
the process coincides with the one described in the beginning (see Section \ref{Sec:setting} for more details).

This defines a Markov chain $\{\sigma(t)\}_{t \in \Z_+}$ on $\Sc_n$,
and we call it a {\bf congestion dynamics}.
It is {\bf reversible} with respect to $\mu_n$,
i.e.,
$\mu_n(\sigma)p(\sigma, \sigma')=\mu_n(\sigma')p(\sigma', \sigma)$
for all $\sigma, \sigma' \in \Sc_n$.
Moreover, it is {\bf irreducible},
i.e.,
for every pair $\sigma_0$ and $\sigma$,
there exists $t \in \Z_+$
such that $\sigma(t)=\sigma$ and $\sigma(0)=\sigma_0$ with a positive probability
if $n<N=\kappa Q$ (cf.\ Section \ref{Sec:setting}).
Further the chain is {\bf aperiodic} since it has a holding probability $p(\sigma, \sigma)>0$
for each $\sigma\in \Sc_n$.
This implies that $\mu_n$ is the unique stationary distribution for the chain,
and the distribution for $\sigma(t)$ with an arbitrary initial condition converges to $\mu_n$
as $t$ tends to infinity.
Let $\Pb_{\sigma_0}$ be the distribution of the process $\{\sigma(t)\}_{t\in \Z_+}$ with the initial state $\sigma_0$.

Let us define the (worst case) total variation {\bf mixing time} for $\{\sigma(t)\}_{t \in \Z_+}$ on $\Sc_n$.
Let
\[
D_\TV(t):=\max_{\sigma_0 \in \Sc_n} D_\TV^{\sigma_0}(t) \quad \text{where $D_\TV^{\sigma_0}(t):=\max_{A \subset \Sc_n}|\Pb_{\sigma_0}(\sigma(t)\in A)-\mu_n(A)|$}.
\]
For each $\e \in (0, 1)$, 
we define
\[
T^\mix(\e):=\min\big\{t \ge 0 \ : \ D_\TV(t)\le \e\big\}.
\]

\begin{theorem}\label{Thm:Gibbs}
Let $Q$ be a fixed integer at least two, $\rho$ be a fixed real in $(0, 1/2]$.
For a positive integer $\kappa$, let $n=\floor{\rho \kappa Q}$.
For the Markov chain $\{\sigma(t)\}_{t\in \Z_+}$ on $\Sc_n$,
it holds that for each $\e \in (0, 1)$,
there exist constant $C_\e$ and $N_\e$ such that for all $n \ge N_\e$,
\[
\left|T^\mix(\e)-\frac{1}{2}n\log n\right|\le C_\e n.
\]
\end{theorem}

This in particular establishes cutoff at time $(1/2)n\log n$ with window of order $n$.
The term which depends on $\e$ is not in the dominating order $n\log n$ but within a smaller order $n$.
For background and other basic results, see \cite[Chapter 18]{LP}.
We discuss more about this phenomenon and related results in Section \ref{Sec:background} below.

Associated to the chain of configurations,
we consider its unlabeled version,
namely,
we only count the number of players who adapt resources.
Let us define the {\bf load profile}
\[
\xi(t):=(\xi_v(t))_{v \in S} \quad \text{where $\xi_v(t):=\xi_v(\sigma(t))$ for each $v \in S$}.
\]
This defines a Markov chain on the set
\[
\Xi_{n, S, \kappa}:=\left\{\xi=(\xi_v)_{v\in S}\in \{0, 1, \dots, \kappa\}^S \ : \ \sum_{v \in S}\xi_v=n\right\}.
\]
The chain is irreducible if $n<N=\kappa Q$,
aperiodic and has a unique stationary distribution $\pi_n$.
Let $\Pb_{\xi_0}$ be the distribution of the process $\{\xi(t)\}_{t\in \Z_+}$ with the initial state $\xi_0$ on $\Xi_{n, S, \kappa}$.
Let us consider the worst case total variation mixing time for $\{\xi(t)\}_{t \in \Z_+}$ against $\pi_n$ on $\Xi_{n, S, \kappa}$ analogously to the case of $\{\sigma(t)\}_{t\in \Z_+}$.

\begin{theorem}\label{Thm:load}
Let $Q$ be a fixed integer at least two,
$\rho$ be a fixed real in $(0, 1/2]$.
For a positive integer $\kappa$,
let $n=\floor{\rho \kappa Q}$.
For the Markov chain $\{\xi(t)\}_{t \in \Z_+}$ on $\Xi_{n, S, \kappa}$,
it holds that for each $\e \in (0, 1)$,
there exist constants $C_\e$ and $N_\e$ such that for all $n \ge N_\e$,
\[
\left|T^\mix(\e)-\frac{1}{2}(1-\rho)n\log n\right| \le C_\e n.
\]
\end{theorem}

The load profile chain exhibits cutoff but at time different from the original congestion dynamics of configurations.
In fact, Theorem \ref{Thm:load} can also be deduced as a time change of the exclusion processes of unlabeled particles on complete graphs \cite[Theorem 1.1]{LacoinLeblond} (see discussions in Section \ref{Sec:background}).
It is illustrative to write the statement in its equivalent form:
For each $\e \in (0, 1)$,
\[
T^\mix(\e)=\frac{1}{2}(1-\rho)\rho N\log N+O_\e(N) \quad \text{as $N \to \infty$}.
\]
This reveals the equivalence between the numbers of players $n$ and $N-n$.
Namely,
considering the process $(\kappa-\xi_v(t))_{v\in S}$ for $t \in \Z_+$ amounts to the load profile chain with $N-n$ players under the same constraint $\kappa$ on each resource.
Thus there is no loss of generality to assume that $\rho \in (0, 1/2]$.
The Glauber dynamics itself, however, does not (at least apparently) have the same kind of equivalence.
We leave open as to whether Theorem \ref{Thm:Gibbs} holds for all $\rho \in (0, 1)$.

\subsection{Background and related results}\label{Sec:background}

This work arises from an attempt to find an explicit class of examples which exhibit cutoff 
in sampling problems of (discrete counterparts of) log-concave distributions and their generalizations.
There has been a significant progress on understanding mixing times of natural Markov chains with stationary distributions whose generating polynomials are (strongly) log-concave \cite{logconcaveII} and \cite{LogconcaveIIannales}.
An important class includes the base-exchange chain of a given matroid $\Mcc=([N], \Ic)$ on the ground set $[N]$ of rank $n$.
It has been shown that for this chain $T^\mix(\e)=O(n\log (n/\e))$ for $\e \in (0, 1)$ \cite[a special case of Theorem 1]{logconcaveIV} developed upon \cite{CryanGuoMousa2021} and \cite{logconcaveII}.
The order is optimal since there is a lower bound $\Omega(n \log n)$, which follows from a coupon collector problem.
Given the results, it is tempting to understand a cutoff for this type of chains,
or to find a reasonably wide class of chains for which we are able to establish an abrupt convergence.

In view of generalizations and algorithmic perspective,
it is natural to consider sampling problems from log M-convex sets or its distribution counterparts in the discrete convex analysis \cite{Murota}.
We say that a probability measure $\pi$ on $\Z^Q$ is {\bf log M-concave} if $\pi$ is supported in the set 
$\{\xi \in \Z^Q \ : \ \sum_{i=1}^Q\xi_i=n\}$ for some $n \in \Z_+$ and the following holds:
For all $\xi$ and $\tilde \xi$ in the support and for every $i \in [Q]$ with $\xi_i>\tilde \xi_i$,
there exists $j \in [Q]$ such that $\xi_j<\tilde \xi_j$ and
\[
\log \pi\big(\xi\big)+\log \pi\big(\tilde \xi\,\big)\le \log \pi\big(\xi-\mb{e}_i+\mb{e}_j\big)+\log \pi\big(\tilde \xi-\mb{e}_j+\mb{e}_i\big).
\]
In the above, $\mb{e}_i$, $i=1, \dots, Q$, denotes the standard basis in $\Z^Q$.
The load profile chain on $\Xi_{n, S, \kappa}$ has the stationary distribution $\pi_n$,
which is obtained by the product of binomial distributions conditioned on $\sum_{v \in S}\xi_v=n$ (cf.\ Section \ref{Sec:setting}).
The measure $\pi_n$ is log M-concave as it is directly checked by this explicit form.
More generally,
if the function $m \mapsto \log (m!)+ \sum_{i=0}^m\f(i)$ is convex,
then $\pi_n$ is log M-concave.
The function $\f$ may also vary in $v \in S$;
see \cite[Example 6.3, Chapter 6]{Murota} for a thorough discussion.
It would be interesting to generalize our result in this setting.
We, however, point out that a cutoff does not necessarily occur for an analogous Markov chain on an M-convex set with the uniform stationary distribution.
This corresponds to the case when $\log \pi_n$ is constant on $\Xi_{n, S, \kappa}$ and $n \le \kappa$ (see Section \ref{Sec:appendix}).

From the viewpoint of interacting particle systems,
the load profile chain is a generalization of exclusion process,
where each site is occupied by at most one indistinguishable particle.
Concerning exclusion processes, mixing times on general underlying graphs (not necessarily complete graphs) have been studied in \cite{Oliveira}.
For relaxation times in this context,  see a recent work \cite{KanegaeWachi}.
In the case of complete graphs on $n$ vertices
(or, equivalently the Bernoulli-Laplace urn model),
the exclusion process of $k$-(unlabeled) particles is defined as follows:
Starting from a configuration of particles where each site has at most one particle on the graph,
independently in each step we choose an (undirected) edge uniformly at random and swap the contents of extremes of the edge.
For the unlabeled particles,
a cutoff has been established by \cite[Theorem 1]{DiaconisShahshahaniBernoulliLaplace}
(for the case of $\floor{n/2}$ particles),
and by \cite[Theorem 1.1]{LacoinLeblond}.
The proof of the former is representation theoretic whereas that of the latter is probabilistic.
In the latter, it has been shown that for the case of $k_n \in [0, n]$ particles and $k_n/\sqrt{n} \to \infty$ as $n \to \infty$,
a cutoff occurs at time $(1/4)n\log n$ with window of order $n$.
Also, in the case of $k_n/\sqrt{n}\to 0$ and $k_n\to \infty$ as $n \to \infty$,
or in the case of $k_n/\sqrt{n} \to l$ for some $l\in (0, \infty)$ as $n \to \infty$,
a cutoff occurs at time $(1/2)n\log n$ with window of order $n$.
Recently, it is announced in \cite{Olesker-Taylor-Schmid}
that the (different) limiting profiles have been identified in these three cases.
See also \cite{NestoridiOleskerTaylor2022} for other results and recent developments in this direction.

The cutoff for the load profile chain in Theorem \ref{Thm:load} actually follows by a time change from the above mentioned result in \cite[Theorem 1.1]{LacoinLeblond}.
Indeed, let us consider the exclusion process of $n$-unlabeled particles on the complete graph on $N$-vertices, and $n/N \to \rho \in (0, 1)$, in particular $n/\sqrt{N} \to \infty$ as $N\to \infty$.
Further each vertex is colored by one of $Q$ different colors and the number of vertices with a common color is $\kappa$ and $N=\kappa Q$.
Independently in each step,
if a chosen edge has the extremes with the same contents (i.e., either both have particles or both have no particles),
then that step is not counted.
This possible skip does not change the
configuration since particles are unlabeled.
Under this time change, the resulting chain is essentially the load profile chain itself.
This speeded-up exclusion process has the mixing time asymptotically equal to
\[
\frac{2n(N-n)}{N(N-1)}\frac{1}{4}N\log N \sim \frac{1}{2}\rho(1-\rho)N\log N.
\]
We produce a self-contained proof of Theorem \ref{Thm:load} adapted to the present setting.
Some of the technical results obtained along the proof are used to show Theorem \ref{Thm:Gibbs}.

The exclusion process with labeled particles has been studied in \cite[Theorem 1.2]{LacoinLeblond} and \cite{PatkoPete}.
A different generalization of the Bernoulli-Laplace urn model (for many urns) is discussed in \cite{Scarabotti1997}.

\subsection{More on the setting}\label{Sec:setting}

Let us describe more in details on the Glauber dynamics in Theorem \ref{Thm:Gibbs}.
The chain is irreducible if $n<\kappa Q$.
Indeed, this follows since given $\sigma \in \Sc_n$ and a distinct pair $i, j \in [n]$
the configuration $\sigma' \in \Sc_n$ obtained from $\sigma$ by exchanging $\sigma_i$ and $\sigma_j$ can be reached within at most $3$ steps.
The transition probability is defined by
\[
p(\sigma, \sigma')=\frac{1}{n}\frac{e^{-\F(\sigma')}}{\sum_{\tau_j=\sigma_j, j \in [n]\setminus\{i\}}e^{-\F(\tau)}} \quad \text{for $\sigma, \sigma' \in \Sc_n$},
\]
if $\sigma_j=\sigma_j'$ for $j \in [n]\setminus\{i\}$,
and $p(\sigma, \sigma')=0$ otherwise.
Recall that $\F(\sigma)=\sum_{v \in S}\sum_{s=0}^{\xi_v(\sigma)}\f(s)$ and 
$\f(s)=\log(\kappa+1)-\log(\kappa+1-s)$.
If $\sigma_i=v_0$ and $\sigma'_i=v_1$ with $v_0\neq v_1$,
then
\[
p(\sigma, \sigma')=\frac{1}{n}\frac{e^{-\f(\xi_{v_1}(\sigma)+1)}}{e^{-\f(\xi_{v_0}(\sigma))}+\sum_{v \neq v_0}e^{-\f(\xi_v(\sigma)+1)}}
=\frac{\kappa-\xi_{v_1}(\sigma)}{n(1+N-n)}.
\]
Furthermore,
we have
\[
p(\sigma, \sigma)=\sum_{i=1}^n\frac{1+\kappa-\xi_{\sigma_i}(\sigma)}{n(1+N-n)}.
\]

For the load profile chain in Theorem \ref{Thm:load},
each transition from $\xi$ to $\xi'$ occurs in the following way:
First $v_-\in S$ is chosen with probability $\xi_{v_-}/n$.
Then we have $\xi'=\xi$ with probability $(1+\kappa-\xi_{v_-})/(1+N-n)$,
or $v_+\in S\setminus\{v_-\}$ is chosen with probability
$(\kappa-\xi_{v_+})/(1+N-n)$ and obtain $\xi'$ from $\xi$ by decrementing $\xi_{v_-}$ by $1$ and then incrementing $\xi_{v_+}$ by $1$.
Equivalently, 
choose $v_-\in S$ with probability $\xi_{v_-}/n$.
We define $\xi'=\xi$ with probability $1/(1+N-n)$ and otherwise define $\xi'$ as follows:
Choose $v_+\in S$ independently with probability $(\kappa-\xi_{v_+})/(N-n)$.
If $v_-=v_+$, then we define $\xi'=\xi$, and if $v_-\neq v_+$,
then
\[
\xi_v':=
\begin{cases}
\xi_v-1 & \text{if $v=v_-$},\\
\xi_v+1 & \text{if $v=v_+$, and}\\
\xi_v & \text{else}.
\end{cases}
\]
Note that the holding (i.e., $\xi'=\xi$) occurs with probability
\[
\sum_{v \in S}\frac{\xi_v}{n}\frac{1+\kappa-\xi_v}{1+N-n}.
\]
Let us define the distribution $\pi_n$ on $\Xi_{n, S, \kappa}$ by
\[
\pi_n(\xi):=\frac{1}{Z_{n, S, \kappa}}\prod_{v \in S}\frac{1}{\xi_v !(\kappa-\xi_v)!} \quad \text{where $Z_{n,S,\kappa}:=\sum_{\xi \in \Xi_{n, S, \kappa}}\prod_{v \in S}\frac{1}{\xi_v!(\kappa-\xi_v)!}$}.
\]
In the above we understand $0!=1$.
The chain $\{\xi(t)\}_{t\in \Z_+}$ is reversible with respect to $\pi_n$.
Indeed, for each distinct pair $\xi$ and $\xi'$ with $\xi_{v_-}'=\xi_{v_-}-1$, 
$\xi_{v_+}'=\xi_{v_+}+1$ and $\xi_v'=\xi_v$ for $v \neq v_-, v_+$,
we have
\[
\frac{\xi_{v_-}}{n}\frac{\kappa-\xi_{v_+}}{1+N-n}{\kappa \choose \xi_{v_-}}{\kappa \choose \xi_{v_+}}
=\frac{\xi_{v_+}+1}{n}\frac{\kappa-\xi_{v_-}+1}{1+N-n}{\kappa \choose \xi_{v_+}+1}{\kappa \choose \xi_{v_-}-1}.
\]
This implies the required reversibility.

\subsection{Outlines of the proofs}

First we discuss Theorem \ref{Thm:load} on the load profile chain and then mention Theorem \ref{Thm:Gibbs} on the congestion dynamics.
The latter uses some of results obtained along the proof of the former.

The proof of Theorem \ref{Thm:load} on the load profile chain proceeds as follows.
Recall that $n=\floor{\rho N}$ and $\rho \in (0, 1/2]$.
For this chain,
the stationary distribution $\pi_n$ is the product of $Q$ copies of binomial distributions on $\{0, \dots, \kappa\}$ conditioned on that $\sum_{v \in S}\xi_v=n$.
Under $\pi_n$,
the normalized vector $(\xi_v/N)_{v \in S}$ is concentrated around the constant vector $(\rho/Q)_{v \in S}$ within the order $O\big(1/\sqrt{N}\big)$.
For the normalized load profile chain,
the time to reach the constant vector within that order is $(1/2)(1-\rho)\rho N\log N$ (Lemma \ref{Lem:square_root}).
After that period,
we construct a coupling to coalesce the chain and another chain starting from the stationary distribution.
The coupling period takes $O(n)$ with high probability (Lemma \ref{Lem:coupling}).
This proves the desired upper bound for the mixing time (Theorem \ref{Thm:load_upper}).
For the lower bound,
we exhibit a worst initial state yielding the tight bound (Theorem \ref{Thm:load_lower}).

The proof of Theorem \ref{Thm:Gibbs} on the congestion dynamics elaborates this analysis of the load profile chain.
We consider the {\it load matrix chain} $(\xi_{u, v}(\sigma(t)))_{u, v \in S}$.
This counts in a configuration $\sigma(t)$ the number of sites in $[n]$ where the state is $u$ at time $0$ and the state is $v$ at time $t$ (see Section \ref{Sec:congestion}).
After a short $O(n)$ time,
it boils down to consider this load matrix chain by exploiting a mean-field nature in this model (Lemma \ref{Lem:TV_stationary}).
Analyzing the dynamics for an appropriately normalized load matrix chain,
we show that the time $(1/2)n\log n$ is needed to reach a typical matrix within order $O\big(1/\sqrt{n}\big)$ (Lemma \ref{Lem:load_matrix_square_root}).
After that period,
we construct a coupling which coalesces two processes in $O(n)$ time with high probability,
combined with the process for the number of vacant sites $(\kappa-\xi_v)_{v \in S}$ (Lemma \ref{Lem:coupling_matrix}).
This coupling is developed upon the one for the load profile chain in Section \ref{Sec:load}.
This leads to the required upper bound for the mixing time (Theorem \ref{Thm:upper}).
For the lower bound, we exhibit a worst initial state and show that the mentioned analysis is optimal by using a variance bound for the load matrix chain (Theorem \ref{Thm:lower}).

\subsection{Organization}

In Section \ref{Sec:load},
we discuss the load profile chain and show Theorem \ref{Thm:load}, based on Theorem \ref{Thm:load_upper} (upper bound) and Theorem \ref{Thm:load_lower} (lower bound).
In Section \ref{Sec:congestion},
we focus on the congestion dynamics and prove Theorem \ref{Thm:Gibbs},
based on Theorem \ref{Thm:upper} (upper bound) and Theorem \ref{Thm:lower} (lower bound).
In Appendix \ref{Sec:appendix}, we introduce a reversible Markov chain on an M-convex set with the uniform stationary distribution, and discuss
an example with no-cutoff.

\subsection*{Notation}
We write $c$, $C$, $C'$, $C''$, \dots, for constants whose exact values possibly change from line to line,
and for $c_\delta$ a constant which depends on the other constant $\delta$ to indicate its dependency.
For real valued sequences $\{f(n)\}_{n \in \Z_+}$ and $\{g(n)\}_{n\in \Z_+}$ such that $g(n)>0$ on non-negative integers $\Z_+$,
we write $f(n)=O(g(n))$
if there exist constants $C\ge 0$ and $N$ such that
$|f(n)| \le C g(n)$ for all $n \ge N$.
Further we write $f(n)=o(g(n))$ if $f(n)/g(n) \to 0$ as $n \to \infty$.
For a real $x$, 
we mean by $\floor{x}$ the largest integer at most $x$ and by $\ceil{x}$ the smallest integer at least $x$.
For a set $A$, we write $A^{\sf c}$ for the complement set and $\#A$ for the cardinality.

\section{The load profile chain}\label{Sec:load}

Let us consider the load profile chain $\{\xi(t)\}_{t\in \Z_+}$ on $\Xi_{n, S, \kappa}$ and recall $N=\kappa Q$.
The associated filtration is denoted by $\{\Fc_t\}_{t\in \Z_+}$, i.e., $\Fc_t$ is the sigma algebra generated by $\xi(s)$ for $s=0, \dots, t$.
For each $v \in S$ and each $t \in \Z_+$,
we have
\[
\xi_v(t+1)-\xi_v(t)=\frac{n-\xi_v(t)}{n}\frac{\kappa-\xi_v(t)}{1+N-n}-\frac{\xi_v(t)}{n}\frac{N-n-(\kappa-\xi_v(t))}{1+N-n}+M_{t+1},
\]
where $|M_{t+1}|\le 2$ and $\Eb[M_{t+1}\mid \Fc_t]=0$ almost surely.
Rearranging the terms,
we obtain
\begin{equation}\label{Eq:load_eq0}
\xi_v(t+1)-\xi_v(t)=\frac{\kappa}{1+N-n}-\frac{N}{1+N-n}\frac{\xi_v(t)}{n}+M_{t+1}.
\end{equation}
For the discussion below, 
it is useful to write the equation in the following way:
For each $v \in S$ and each $t \in \Z_+$,
\[
\xi_v(t+1)-\xi_v(t)=\frac{\kappa}{N-n}-\frac{N}{N-n}\frac{\xi_v(t)}{n}+\wh M_{t+1},
\]
where $|\wh M_{t+1}|\le C$ and $|\Eb[\wh M_{t+1}\mid \Fc_t]|\le C/N$ almost surely.
Further $C$ is a constant depending only on $\rho$.
Note that $\rho-1/N\le n/N \le \rho$ and that $N=\kappa Q$ where $Q$ is a fixed constant.
Abusing the notation $\wh M_{t+1}$ and changing the constant $C$,
we have
\[
\xi_v(t+1)-\xi_v(t)=\frac{1}{(1-\rho)Q}-\frac{1}{(1-\rho)\rho N}\xi_v(t)+\wh M_{t+1}.
\]
This leads to
\begin{equation}\label{Eq:load_eq}
\frac{\xi_v(t+1)}{N}-\frac{\rho}{Q}=\(1-\frac{1}{(1-\rho)\rho N}\)\(\frac{\xi_v(t)}{N}-\frac{\rho}{Q}\)+\frac{\wh M_{t+1}}{N}.
\end{equation}

Simplifying notations,
we write
\[
\Delta(t):=(\Delta_v(t))_{v\in S} 
\quad
\text{where $\Delta_v(t):=\frac{\xi_v(t)}{N}-\frac{\rho}{Q}$}.
\]
Let us define
\[
\gamma:=\frac{1}{(1-\rho)\rho}.
\]
Squaring both sides of \eqref{Eq:load_eq} 
and summing over $v \in S$ yields the following:
\begin{equation*}
\|\Delta(t+1)\|^2=\(1-\frac{\gamma}{N}\)^2\|\Delta(t)\|^2+M_{t+1}',
\end{equation*}
where $|M_{t+1}'|\le CQ/N$ and $|\Eb[ M_{t+1}'\mid \Fc_t]|\le CQ/N^2$ almost surely.
In the above $\|\cdot\|$ denotes the standard Euclidean norm  in $\R^Q$.
Rearranging terms leads to
\begin{equation}\label{Eq:load}
\|\Delta(t+1)\|^2=\(1-\frac{2\gamma}{N}\)\|\Delta(t)\|^2+\wt M_{t+1}.
\end{equation}
In the above $|\wt M_{t+1}|\le CQ/N$ and $|\Eb[\wt M_{t+1}\mid \Fc_t]|\le CQ/N^2$ almost surely.
This shows that for each $t \in \Z_+$, almost surely
\begin{equation}\label{Eq:load_cond}
\Eb[\|\Delta(t+1)\|^2 \mid \Fc_t]\le \(1-\frac{2\gamma}{N}\)\|\Delta(t)\|^2+\frac{CQ}{N^2}.
\end{equation}

The following is a general result which we use throughout.

\begin{lemma}\label{Lem:general}
Let $\{x_t\}_{t\in \Z_+}$
be real valued random variables with the associated filtration $\{\Fc_t\}_{t \in \Z_+}$.
Suppose that there exist constants $\lambda, C>0$ such that for all $t \in \Z_+$ and $N>0$,
\[
x_{t+1}\le \(1-\frac{\lambda}{N}\)x_t+M_{t+1}
\quad
\text{where $|M_{t+1}|\le \frac{C}{N}$ and $|\Eb[M_{t+1}\mid \Fc_t]| \le \frac{C}{N^2}$ almost surely}.
\]
Then there exist constants $c, R_0>0$ such that
for all $R\ge R_0$, $t \in \Z_+$ and $N>0$,
\[
\Pb_{x_0}\(x_t \ge \(1-\frac{\lambda}{N}\)^t x_0+\frac{R}{\sqrt{N}}\)\le e^{-cR^2}.
\]
\end{lemma}

\proof
Let us define inductively $\wbar x_0:=x_0$ and for $t\in \Z_+$,
\begin{equation}\label{Eq:Lem:general}
\wbar x_{t+1}:=\(1-\frac{\lambda}{N}\)\wbar x_t+M_{t+1}.
\end{equation}
It holds that $x_t\le \wbar x_t$ almost surely for each $t \in \Z_+$.
Hence it suffices to show the claim for $\{\wbar x_t\}_{t\in \Z_+}$, and thus we assume that $\{x_t\}_{t\in \Z_+}$ satisfies \eqref{Eq:Lem:general}.
By assumption,
almost surely
\[
\Eb[x_{t+1}\mid \Fc_t]\le \(1-\frac{\lambda}{N}\)x_t+\frac{C}{N^2}.
\]
Therefore it holds that almost surely
\begin{align*}
\Eb\left[x_{t+1}-\frac{C}{N\lambda} \,\Big|\, \Fc_t\right]
&\le \(1-\frac{\lambda}{N}\)x_t-\frac{C}{N \lambda}+\frac{C}{N^2}\\
&=\(1-\frac{\lambda}{N}\)x_t-\frac{C(1-\lambda/N)}{N\lambda}
=\(1-\frac{\lambda}{N}\)\(x_t-\frac{C}{N\lambda}\).
\end{align*}
Letting $z_t:=(1-\lambda/N)^{-t}(x_t-C/(N\lambda))$,
we find that $\{z_t\}_{t\in \Z_+}$ forms a supermartingale with respect to $\{\Fc_t\}_{t\in \Z_+}$.
Note that by \eqref{Eq:Lem:general},
\begin{align*}
|z_{t+1}-z_t|
&=\left|\(1-\frac{\lambda}{N}\)^{-t-1}M_{t+1}-\(1-\frac{\lambda}{N}\)^{-t-1}\frac{C}{N\lambda}+\(1-\frac{\lambda}{N}\)^{-t}\frac{C}{N\lambda}\right|\\
&=\(1-\frac{\lambda}{N}\)^{-t-1}\left|M_{t+1}-\frac{C}{N^2}\right|
\le \(1-\frac{\lambda}{N}\)^{-t-1}\frac{2C}{N}.
\end{align*}
Thus the Azuma-Hoeffding inequality yields for all $R \ge 0$,
\begin{align*}
&\Pb_{x_0}\(z_t \ge z_0+\(1-\frac{\lambda}{N}\)^{-t}R\)
\le \exp\(-\frac{(1-\lambda/N)^{-2t}R^2}{2\sum_{s=0}^t(1-\lambda/N)^{-2s-2}(2C/N)^2}\)\\
&\qquad \qquad \qquad \qquad\le \exp\(-\frac{R^2N^2}{8C^2}(1-\lambda/N)^4((1-\lambda/N)^{-2}-1)\)\le \exp(-cR^2N).
\end{align*}
Replacing $R$ by $R/\sqrt{N}$ in the above inequalities, 
we obtain
\[
\Pb_{x_0}\(x_t-\frac{C}{N\lambda} \ge \(1-\frac{\lambda}{N}\)^t\(x_0-\frac{C}{N\lambda}\)+\frac{R}{\sqrt{N}}\) \le e^{-cR^2}.
\]
Shifting $R$ to $R-C/(\lambda \sqrt{N})$ and choosing a smaller $c>0$ yields the required inequality.
\qed

\medskip

For $\delta>0$,
let
\[
\Xi(\delta):=\left\{\xi \in \Xi_{n, S, \kappa} \ : \ \left\|\frac{\xi}{N}-\(\frac{\rho}{Q}\)_{v \in S}\right\|\le \delta\right\}.
\]

\begin{lemma}\label{Lem:load_good}
For every $\delta>0$,
there exist constants $c_\delta>0$ and $N_\delta$ such that for all $\xi_0 \in \Xi(\delta)$
and for all $N \ge N_\delta$,
\[
\Pb_{\xi_0}\Bigg(\bigcup_{0\le t \le N^2}\Big\{\xi(t) \notin \Xi(2\delta)\Big\}\Bigg) \le e^{-c_\delta N}.
\]
\end{lemma}

\proof
Applying to Lemma \ref{Lem:general} by setting $x_t=\|\Delta(t)\|^2$ and $\lambda=2\gamma$ in \eqref{Eq:load} shows that there exist constants $c, R_0>0$ such that for all $R \ge R_0$ and $N>0$,
\[
\Pb_{\xi_0}\(\|\Delta(t)\|^2\ge \(1-\frac{2\gamma}{N}\)^t\|\Delta(0)\|^2+\frac{R}{\sqrt{N}}\) \le e^{-cR^2}.
\]
Thus for $\delta>0$,
letting $R:=\delta^2\sqrt{N}$ and $N_\delta:=\ceil{(R_0/\delta^2)^2}$,
we find that if $\|\Delta(0)\| \le \delta$,
then for all $N \ge N_\delta$,
\[
\Pb_{\xi_0}\(\|\Delta(t)\|^2 \ge 2 \delta^2\)\le e^{-c\delta^4 N}.
\]
By the union bound over $0\le t \le N^2$,
we have a constant $c_\delta>0$ such that for all $N \ge N_\delta$,
\[
\Pb_{\xi_0}\(\max_{0\le t\le N^2}\|\Delta(t)\|^2 \ge 2\delta^2\) \le e^{-c_\delta N}.
\]
This implies the claim.
\qed

\begin{lemma}\label{Lem:square_root}
Let $t_N:=\ceil{(1/2)(1-\rho)\rho N \log N}$.
There exists a constant $C$ such that for all $N>0$, all real $R>0$ and all $\xi_0 \in \Xi_{n, S, \kappa}$,
\[
\Pb_{\xi_0}\(\xi(t_N) \notin \Xi\(\frac{R}{\sqrt{N}}\)\) \le \frac{CQ}{R^2}.
\]
\end{lemma}

\proof
Taking expectations on both sides of \eqref{Eq:load_cond}
provides the following by induction in $t \in \Z_+$:
There exists $C>0$ such that for all $N>0$ and all $t \in \Z_+$,
\[
\Eb\|\Delta(t)\|^2\le \(1-\frac{2\gamma}{N}\)^t\Eb\|\Delta(0)\|^2+\frac{CQ}{2\gamma N}.
\]
Since $\|\Delta(0)\|^2 \le Q$ and $t_N=\ceil{(2\gamma)^{-1}N\log N}$,
this implies that for $N>0$,
\[
\Eb\|\Delta(t_N)\|^2 \le Q\(1-\frac{2\gamma}{N}\)^{(2\gamma)^{-1}N\log N}+\frac{CQ}{2\gamma N}
\le \frac{Q}{N}+\frac{CQ}{2\gamma N}=\frac{Q}{N}\(1+\frac{C}{2\gamma}\).
\]
Hence the Chebyshev inequality yields
\[
\Pb_{\xi_0}\(\|\Delta(t_N)\| \ge \frac{R}{\sqrt{N}}\) \le \frac{N}{R^2}\Eb\|\Delta(t_N)\|^2\le \frac{Q}{R^2}\(1+\frac{C}{2\gamma}\),
\]
implying the claim.
\qed

\medskip

For $\xi_0, \tilde \xi_0 \in \Xi_{n, S, \kappa}$,
let us consider chains $\{\xi(t)\}_{t\in \Z_+}$ starting from $\xi_0$ and $\{\tilde \xi(t)\}_{t\in \Z_+}$ starting from $\tilde \xi_0$.
We construct a coupling between them such that their coalescent time is at most $O(N)$ with high probability 
if $\xi_0, \tilde \xi_0 \in \Xi(R/\sqrt{N})$ for $R>0$.
Let
\[
D(\xi, \tilde \xi):=\frac{1}{2}\sum_{v \in S}|\xi_v-\tilde \xi_v| \quad \text{for $\xi, \tilde \xi \in \Xi_{n, S, \kappa}$}.
\]
For $t\in \Z_+$, we write $D_t:=D\big(\xi(t), \tilde \xi(t)\big)$ for brevity.

\begin{lemma}\label{Lem:coupling}
Let us fix a constant $0<\delta<\min\{\rho, 1-\rho\}/(4Q)$.
There exists a coupling $\Pb$ between $\{\xi(t)\}_{t\in \Z_+}$ and $\{\tilde \xi(t)\}_{t\in \Z_+}$ satisfying the following:
For each $t\in \Z_+$,
if $\xi(t), \tilde \xi(t) \in \Xi(2\delta)$ and $D_t>0$,
then almost surely for each $r=1$, $-1$,
\[
\Pb\(D_{t+1}-D_t=r \,\Big |\, \xi(t), \tilde \xi(t)\)\ge \frac{1}{4Q^2}
\quad \text{and} \quad
\Eb\left[D_{t+1}-D_t \, \Big|\, \xi(t), \tilde \xi(t)\right] \le 0.
\]
\end{lemma}

\proof
Given $\xi, \tilde \xi \in \Xi_{n, S, \kappa}$,
if $\xi, \tilde \xi \in \Xi(2\delta)$,
then we run the pair of chains $\{(\xi(t), \tilde \xi(t))\}_{t \in \Z_+}$ by the transition defined below,
otherwise we run them independently.
Let $\xi, \tilde \xi \in \Xi(2\delta)$.
If $D(\xi, \tilde \xi)=0$,
then we define common transitions for $\xi$ and $\tilde \xi$,
and $D(\xi, \tilde \xi)$ remains $0$.
Let us assume that $D(\xi, \tilde \xi)>0$.
Then, 
there exists a pair $(v_0, v_1) \in S\times S$ such that
\[
\xi_{v_0}<\tilde \xi_{v_0} \quad \text{and} \quad \xi_{v_1}>\tilde \xi_{v_1}.
\]
Let us fix such a pair $(v_0, v_1)$.
For $v \in S$,
let
\[
m_v:=\min\{\xi_v, \tilde \xi_v\} \quad \text{and} \quad M_v:=\max\{\xi_v, \tilde \xi_v\}.
\]
Since $\xi, \tilde \xi \in \Xi(2\delta)$,
it holds that for all $v \in S$,
\[
m_v \ge N\(\frac{\rho}{Q}-2\delta\)\ge \frac{\rho N}{2Q}
\quad \text{and} \quad
\kappa-M_v=\min\{\kappa-\xi_v, \kappa-\tilde \xi_v\}\ge \frac{(1-\rho)N}{2Q}.
\]

\begin{figure}
\centering
\begin{tikzpicture}[scale=0.5]


\draw (1,6.5) node {{$O_v^o$}};
\draw (2.5,6.5) node {{$O_v^r$}};
\draw (4,6.5) node {{$O_v^+$}};
\draw (5.5,6.5) node {{$V_v^r$}};
\draw (7,6.5) node {{$V_v^o$}};

\draw (12,6.5) node {{$V_w^+$}};
\draw (12,-0.7) node {{$\wt O_w^+$}};


\fill [fill=gray] (0,0)--(0,6)--(2,6)--(2,0)--cycle;
\fill [fill=gray] (3,0)--(3,6)--(5,6)--cycle;
\fill [fill=lightgray] (3,0)--(5,0)--(5,6)--cycle;
\fill [fill=lightgray] (6,0)--(6,6)--(8,6)--(8,0)--cycle;

\fill [fill=gray] (8,0)--(8,6)--(10,6)--(10,0)--cycle;
\fill [fill=gray] (11,0)--(13,0)--(11,6)--cycle;
\fill [fill=lightgray] (11,6)--(13,6)--(13,0)--cycle;
\fill [fill=lightgray] (14,0)--(14,6)--(16,6)--(16,0)--cycle;

\fill [fill=gray] (16,0)--(16,6)--(18,6)--(18,0)--cycle;
\fill [fill=gray] (19,0)--(21,0)--(19,6)--cycle;
\fill [fill=lightgray] (19,6)--(21,6)--(21,0)--cycle;
\fill [fill=lightgray] (22,0)--(22,6)--(24,6)--(24,0)--cycle;

\draw[very thick] (0,0)--(0,6);

\draw[thick, densely dotted] (2,0)--(2,6);
\draw[solid] (3,0)--(3,6);
\draw[solid] (3,0)--(5,6);
\draw[solid] (5,0)--(5,6);
\draw[thick, densely dotted] (6,0)--(6,6);

\draw[very thick] (8,0)--(8,6);

\draw[thick, densely dotted] (10,0)--(10,6);
\draw[solid] (11,0)--(11,6);
\draw[solid] (13,0)--(11,6);
\draw[solid] (13,0)--(13,6);
\draw[thick, densely dotted] (14,0)--(14,6);

\draw[very thick] (16,0)--(16,6);

\draw[thick, densely dotted] (18,0)--(18,6);
\draw[solid] (19,0)--(19,6);
\draw[solid] (21,0)--(19,6);
\draw[solid] (21,0)--(21,6);
\draw[thick, densely dotted] (22,0)--(22,6);

\draw[very thick] (24,0)--(24,6);

\draw[solid] (0,6)--(24,6);
\draw[solid] (0,0)--(24,0);

\draw (1.7,-0.7) node {{$m_v$}};
\draw (4,-0.7) node {{$\widetilde V_v^+$}};
\draw (6.5,-0.7) node {{$\kappa-M_v$}};

\draw[ultra thick] (0,0)--(3,0);
\draw[ultra thick] (5,0)--(8,0);
\end{tikzpicture}
\caption{An illustration of a sequence $(\mb{s}_i)_{i \in [N]}$ of elements in $(S\sqcup \wbar S)^2$,
where the first coordinates are aligned on the upper row and the second coordinates are aligned on the lower row for $\mb{s}_i=(v_i, v'_i)$ and $i \in [N]$.}
\label{Fig:coupling_load}
\end{figure}
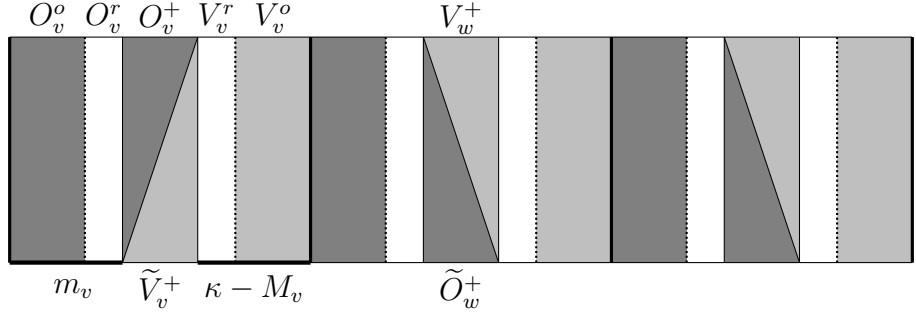

Let us define a sequence $(\mb{s}_i)_{i \in [N]}$ of elements in $(S\sqcup \wbar S)^2$
where $\wbar S:=\{\wbar v \ : \ v \in S\}$.
In this notation,
$v$ and $\wbar v$ indicate ``occupancy'' and ``vacancy'' respectively.
We align $\mb{s}_i \in (S\sqcup\wbar S)^2$ for $i \in [N]$
by counting $(\xi_v)_{v\in S}$ and $(\tilde \xi_v)_{v\in S}$.
See also Figure \ref{Fig:coupling_load}.
For $v \in S=\{1, \dots, Q\}$,
let us define $(\mb{s}_i)_{i \in [N]}$ as in the following:
\[
\mb{s}_i=
\begin{cases}
(v, v) & \text{for $i \in O_v:=[(v-1)\kappa+1, (v-1)\kappa+m_v]$},\\
(v, \wbar v) & \text{for $i \in O_v^+:=[(v-1)\kappa+m_v+1, (v-1)\kappa+\xi_v]$},\\
(\wbar v, v) & \text{for $i \in V_v^+:=[(v-1)\kappa+m_v+1, (v-1)\kappa+\tilde \xi_v]$},\\
(\wbar v, \wbar v) & \text{for $i \in V_v:=[v \kappa+1-(\kappa-M_v), v \kappa]$}.
\end{cases}
\]
Note that either $O_v^+$ or $V_v^+$ is empty for each $v \in S$.
Let us also write $\wt O_v^+:=V_v^+$
and $\wt V_v^+:=O_v^+$.
We further define
\[
O_v^o:=\left[(v-1)\kappa+1, (v-1)\kappa +\ceil{\frac{n}{2Q}}\right]
\quad \text{and} \quad
V_v^o:=\left[v\kappa+1-\ceil{\frac{N-n}{2Q}}, v \kappa\right],
\]
also $O_v^r:=O_v\setminus O_v^o$ and $V_v^r:=V_v\setminus V_v^o$.
These are well-defined when $\xi, \tilde \xi \in \Xi(2\delta)$.
Let us define $O:=\bigsqcup_{v \in S}O_v$ and similarly $V$, $O^o$, $V^o$, $O^r$, $V^r$, $O^+$, $V^+$, $\wt O^+$ and $\wt V^+$ by taking the unions over $v \in S$ in the corresponding cases.
Note that the cardinalities of $O^+$, $V^+$, $\wt O^+$ and $\wt V^+$ coincide and equal $D(\xi, \tilde \xi)$.
Using this fact,
we choose bijections
\[
\theta^O: O^+\to \wt O^+ \quad \text{and} \quad \theta^V:V^+ \to \wt V^+.
\]

We retain the current state with probability $1/(1+N-n)$,
or perform the following with probability $(N-n)/(1+N-n)$.
Choose 
\[
a \in O\sqcup O^+ \quad \text{and} \quad b \in V\sqcup V^+
\]
uniformly at random respectively and independently.
For the sake of brevity,
we say that {\bf $\xi'$ is obtained by applying $(v, w) \in S \times S$ to $\xi$} 
when we define $\xi'$ from $\xi$,
decrementing $\xi_v$ by $1$ and then incrementing $\xi_w$ by $1$.
Note that if $v=w$,
then $\xi'=\xi$.
We define and analyze the transition from $(\xi, \tilde \xi)$ to $(\xi', \tilde \xi')$, depending on the following cases:
\subsubsection*{\upshape{(i) {\bf Case} $(a, b) \in O \times V$}}
If $(a, b) \in O_v \times V_w\setminus O_v^o\times V_w^o$,
then we define $\xi'$ and $\tilde \xi'$ by applying $(v, w)$ to $\xi$ and to $\tilde \xi$ respectively.
Let us consider the case when $(a, b)\in O^o_v\times V^o_w$.
Since the cardinalities are equal, there exist bijections
\[
O_{v_0}^o \times V_{v_0}^o \to O_{v_0}^o\times V_{v_1}^o 
\quad \text{and} \quad
O_{v_1}^o\times V_{v_1}^o \to O_{v_1}^o\times V_{v_0}^o.
\]
For $k=0, 1$,
if $(a, b) \in O_{v_k}^o\times V_{v_k}^o$,
then we define $\xi':=\xi$
and $\tilde \xi'$ by applying $(v_k, v_{1-k})$ to $\tilde \xi$,
and if $(a, b) \in O_{v_k}^o\times V_{v_{1-k}}^o$,
then we define $\xi'$ by applying $(v_k, v_{1-k})$ to $\xi$
and $\tilde \xi'=\tilde \xi$.
If $(a, b) \in O_v^o\times V_w^o$ with $(v, w) \notin \{v_0, v_1\}\times \{v_0, v_1\}$,
then we define $\xi'$ and $\tilde \xi'$ by applying $(v, w)$ to $\xi$ and to $\tilde \xi$ respectively.

In this case $D_{t+1}-D_t=-1$ or $1$,
which happens with equal probability
\[
2\frac{1}{n}\ceil{\frac{\rho N}{2Q}}\frac{1}{N-n}\ceil{\frac{(1-\rho)N}{2Q}}\frac{N-n}{1+N-n}\ge \frac{1}{4Q^2},
\]
and $D_{t+1}-D_t=0$ otherwise since the identical move is applied to both chains.

\subsubsection*{\upshape{(ii) {\bf Case} $(a, b) \in O\times V^+\sqcup O^+\times V$}}
Note that there exist bijections
\[
O\times V^+\to O\times \wt V^+, (a, b)\mapsto (a, \theta^V(b))
\quad \text{and} \quad 
O^+\times V \to \wt O^+\times V,
(a, b)\mapsto (\theta^O(a), b).
\]
If $(a, b)\in O_v\times V_w^+$,
then $(a, \theta^V(b)) \in O_v\times \wt V_{\wt w}^+$ for $\wt w \in S$ and $\wt w\neq w$.
In this case, $\xi_w<\tilde \xi_w$ and $\xi_{\wt w}>\tilde \xi_{\wt w}$.
We define $\xi'$ by applying $(v, w)$ to $\xi$,
and $\tilde \xi'$ by applying $(v, \wt w)$ to $\tilde \xi$.
Similarly, if $(a, b) \in O_v^+\times V_w$,
then $(\theta^O(a), b)\in \wt O_{\wt v}^+\times V_w$ for $\wt v \in S$ and $\wt v\neq v$.
In this case, $\xi_v>\tilde \xi_v$ and $\xi_{\wt v}<\tilde \xi_{\wt v}$.
We define $\xi'$ by applying $(v, w)$ to $\xi$,
and $\tilde \xi'$ by applying $(\wt v, w)$ to $\tilde \xi$.

In this case
$D_{t+1}-D_t=-2$,
which follows since $V_w^+$ and $\wt V_{\wt w}^+$ (respectively, $O_v^+$ and $\wt O_{\wt v}^+$) are non-empty (if this case happens) and $w \neq \wt w$ (respectively, $v \neq \wt v$).
(Note that it is possible that $v=w$ or $v=\wt w$ but not both since $w \neq \wt w$,
similarly it is possible that $v=w$ or $\wt v=w$ but not both since $v\neq \wt v$.)

\subsubsection*{\upshape{(iii) {\bf Case} $(a, b) \in O^+\times V^+$}}
Recall that $O_v^+=\wt V_v^+$ and $V_w^+=\wt O_w^+$ for all $v, w \in S$.
If $(a, b)\in O_v^+\times V_w^+$,
then we define $\xi'$ by applying $(v, w)$ to $\xi$,
and $\tilde \xi'$ by applying $(w, v)$ to $\tilde \xi$.

In this case $D_{t+1}-D_t=-2, -1$ or $0$ 
since $v\neq w$, and the roles of $v$ and $w$ are interchanged for $\xi(t)$ and $\tilde \xi(t)$.
The difference is $-2$ when both of $\#O_v^+$ and $\#V_w^+$ are at least $2$,
it is $-1$ when one of them is $1$ and the other is at least $2$,
and it is $0$ when both of them are $1$.

Summarizing the above cases, we conclude the claim as required.
\qed

\begin{theorem}\label{Thm:load_upper}
Let us fix an integer $Q \ge 2$ and a real $\rho \in (0, 1)$.
For a positive integer $\kappa$,
let $N=\kappa Q$ and $n=\floor{\rho N}$.
For the load profile chain $\{\xi(t)\}_{t\in \Z_+}$ on $\Xi_{n, S, \kappa}$,
letting $t_N:=\ceil{(1/2)(1-\rho)\rho N\log N}$,
we have
\[
\lim_{\alpha\to \infty}\limsup_{N \to \infty}D_\TV\(t_N+\ceil{\alpha N}\)=0.
\]
\end{theorem}

\proof
For each $\xi(0) \in \Xi_{n, S, \kappa}$,
let us define a coupling between the load profile chains $\{\xi(t)\}_{t\in \Z_+}$ and $\{\tilde \xi(t)\}_{t\in \Z_+}$ starting from $\xi(0)$ and from the stationary distribution $\pi_n$ respectively.
We run two chains independently before time $t_N$,
and then apply to the pair the coupling constructed in Lemma \ref{Lem:coupling}.
Let
\[
\tau_{\coup}:=\inf\left\{t \ge 0 \ : \ \xi(t)=\tilde \xi(t)\right\}.
\]
We will estimate the probability that $\tau_{\coup}>t_N+\alpha N$ under the coupling.

For each $R>0$,
let us define the event
\[
A_{N, R}:=\left\{\xi(t_N), \tilde \xi(t_N) \in \Xi\(\frac{R}{\sqrt{N}}\)\right\}.
\]
By Lemma \ref{Lem:square_root},
there exists $C$ such that for all $N>0$ and $R>0$,
\begin{equation}\label{Eq:PAc}
\Pb_{\xi(0), \pi_n}\(A_{N, R}^{\sf c}\)\le \Pb_{\xi(0)}\(\xi(t_N) \notin \Xi\(\frac{R}{\sqrt{N}}\)\)+\Pb_{\pi_n}\(\xi(t_N) \notin \Xi\(\frac{R}{\sqrt{N}}\)\) \le \frac{2CQ}{R^2}.
\end{equation}
Recall that $A^{\sf c}$ of an event $A$ stands for the complement event.

Let us fix a constant $\delta$ such that $0<\delta\le \min\{\rho, 1-\rho\}/(4Q)$.
For each $R>0$, for all large enough $N$,
we have $R/\sqrt{N} \le \delta$.
For $t \in \Z_+$,
let
\[
G_t:=\Big\{\xi(t) \in \Xi(2\delta)\Big\} \quad \text{and} \quad \wt G_t:=\Big\{\tilde \xi(t) \in \Xi(2\delta)\Big\}.
\]
Lemma \ref{Lem:load_good} implies that 
if $\xi_0, \tilde \xi_0 \in \Xi(\delta)$,
for all large enough $N$,
\begin{equation}\label{Eq:G}
\Pb_{\xi_0}\Big(\bigcup_{0\le t \le N^2}G_t^{\sf c}\Big) \le e^{-c_\delta N}
\quad \text{and} \quad 
\Pb_{\tilde \xi_0}\Big(\bigcup_{0\le t \le N^2}\wt G_t^{\sf c}\Big) \le e^{-c_\delta N}.
\end{equation}

By Lemma \ref{Lem:coupling},
there exists a coupling satisfying the following:
If $\xi(t), \tilde \xi(t)\in \Xi(2\delta)$ and $D_t>0$,
then the variance of $D_{t+1}-D_t$ conditioned on $\xi(t)$ and $\tilde \xi(t)$ satisfies that
\[
\Var(D_{t+1}-D_t\mid \xi(t), \tilde \xi(t)) \ge \frac{2}{4Q^2}=\frac{1}{2Q^2}.
\]
Thus the probability that the first hitting time of $\{D_t\}_{t\in \Z_+}$ at $0$ is more than $\alpha N$ for $\alpha>0$ is at most
\[
\frac{4D_0}{\sqrt{\alpha N/(2Q^2)}}\le \frac{6Q D_0}{\sqrt{\alpha N}}.
\]
This follows from \cite[Proposition 17.20]{LP}.
Hence for $\alpha>0$ and $R>0$ for all large enough $N$
and for all $\xi_0, \tilde \xi_0 \in \Xi(R/\sqrt{N})$,
the following holds by \eqref{Eq:G}:
\begin{align*}
\Pb_{\xi_0, \tilde \xi_0}\(\tau_{\coup}>\alpha N\)
&\le \Pb_{\xi_0, \tilde \xi_0}\Big(\{\tau_{\coup}>\alpha N\}\cap \bigcap_{0\le t \le \ceil{\alpha N}}G_t\cap \wt G_t\Big)+2e^{-c_\delta N}\\
&\le \frac{6QD_0}{\sqrt{\alpha N}}+2e^{-c_\delta N}
\le \frac{6Q^2 R}{\sqrt{\alpha}}+2e^{-c_\delta N} \le \frac{10Q^2R}{\sqrt{\alpha}}.
\end{align*}
In the third inequality above we have used $D_0 \le QR\sqrt{N}$, which holds if $\xi_0, \tilde \xi_0 \in \Xi(R/\sqrt{N})$.
This implies that by the Markov property,
\begin{equation}\label{Eq:coupling}
\Pb\(\tau_\coup > t_N+\alpha N \mid A_{N, R}\) \le \frac{10Q^2 R}{\sqrt{\alpha}}.
\end{equation}
Therefore by \eqref{Eq:PAc} and \eqref{Eq:coupling},
we obtain
\begin{align*}
\Pb_{\xi(0), \pi_n}\(\tau_\coup>t_N+\alpha N\)
&\le \Pb_{\xi(0), \pi_n}\(\{\tau_\coup>t_N+\alpha N\}\cap A_{N, R}\)+\frac{2CQ}{R^2}\\
&\le \Pb\(\tau_\coup >t_N+\alpha N \mid A_{N, R}\)+\frac{2CQ}{R^2}
\le \frac{10Q^2R}{\sqrt{\alpha}}+\frac{2CQ}{R^2}.
\end{align*}
Note that the above estimate is uniform over $\xi(0)\in \Xi_{n,S,\kappa}$,
and
\[
\max_{A \subset \Xi_{n, S, \kappa}}\left|\Pb_{\xi(0)}\(\xi(t_N+\ceil{\alpha N})\in A\)-\pi_n\(A\)\right|
\le \Pb_{\xi(0), \pi_n}\(\tau_\coup>t_N+\alpha N\).
\]
For all $\alpha>0$ and $R>0$, and for all large enough $N$,
we have
\[
D_\TV\(t_N+\ceil{\alpha N}\)\le \max_{\xi(0) \in \Xi_{n, S, \kappa}}\Pb_{\xi(0), \pi_n}\(\tau_\coup>t_N+\alpha N\)
\le \frac{10 Q^2 R}{\sqrt{\alpha}}+\frac{2CQ}{R^2}.
\]
Hence, letting $N \to \infty$ and then $\alpha \to \infty$ yields
\[
\limsup_{\alpha \to \infty}\limsup_{N\to \infty}D_\TV\(t_N+\ceil{\alpha N}\) \le \frac{2CQ}{R^2}.
\]
Since the above estimate holds for all $R>0$, we conclude the claim.
\qed

\begin{theorem}\label{Thm:load_lower}
In the same setting as in Theorem \ref{Thm:load_upper},
for the load profile chain $\{\xi(t)\}_{t\in \Z_+}$ on $\Xi_{n, S, \kappa}$,
letting $t_N:=\ceil{(1/2)(1-\rho)\rho N \log N}$,
we have
\[
\lim_{\alpha \to \infty}\liminf_{N \to \infty} D_\TV\(t_N-\ceil{\alpha N}\)=1.
\]
\end{theorem}

\proof
For each $t \in \Z_+$ and each $v \in S$,
by \eqref{Eq:load_eq0},
\[
\frac{\xi_v(t+1)}{N}-\Eb \frac{\xi_v(t+1)}{N}=\(1-\frac{N}{(1+N-n)n}\)\(\frac{\xi_v(t)}{N}-\Eb \frac{\xi_v(t)}{N}\)+\frac{M_{t+1}}{N},
\]
where $|M_{t+1}|\le 2$ and $\Eb[M_{t+1}\mid \Fc_t]=0$ almost surely.
Squaring both sides and taking the expectations leads to
\[
\Var\(\frac{\xi_v(t+1)}{N}\)= \(1-\frac{N}{(1+N-n)n}\)^2\Var\(\frac{\xi_v(t)}{N}\)+\frac{\Eb|M_{t+1}|^2}{N^2}.
\]
By induction in $t \in \Z_+$,
this implies that for a constant $C$ independent of $t$ or $N$,
\[
\Var\(\frac{\xi_v(t)}{N}\) \le \(1-\frac{\gamma}{N}\)^{2t}\Var\(\frac{\xi_v(0)}{N}\)+\frac{C}{(1-(1-\gamma/N)^2)N^2}.
\]
For every fixed initial state $\xi(0) \in \Xi_{n, S, \kappa}$,
we have $\Var(\xi_v(0)/N)=0$,
and thus
\begin{equation}\label{Eq:load_var}
\Var\(\frac{\xi_v(t)}{N}\) \le \frac{C}{(2\gamma/N-\gamma^2/N^2)N^2}\le \frac{C}{\gamma N},
\end{equation}
for all $N \ge \gamma$ and for all $t \in \Z_+$.

Let us fix $\xi(0)$ such that 
$\xi_v(0)=\min\{\kappa, n\}$ for some fixed $v \in S$.
Since $\rho$ and $Q$ are fixed constants,
there exists $c_0>0$ such that for all large enough $N$,
\[
\frac{\xi_v(0)}{N}-\frac{\rho}{Q} \ge c_0.
\]
By \eqref{Eq:load_eq} and by taking expectations, 
it holds that for an independent constant $C>0$,
\[
\Eb \frac{\xi_v(t+1)}{N}-\frac{\rho}{N} \ge \(1-\frac{\gamma}{N}\)\(\Eb\frac{\xi_v(t)}{N}-\frac{\rho}{N}\)-\frac{C}{N^2}.
\]
By induction in $t \in \Z_+$,
we obtain
\[
\Eb \frac{\xi_v(t)}{N}\ge \frac{\rho}{Q}+\(1-\frac{\gamma}{N}\)^t\Eb\(\frac{\xi_v(0)}{N}-\frac{\rho}{Q}\)-\frac{C}{\gamma N}
\ge \frac{\rho}{Q}+\(1-\frac{\gamma}{N}\)^t c_0-\frac{C}{\gamma N}.
\]
Therefore for $T_{N,\alpha}:=\floor{t_N-(\alpha/\gamma)N}$ for $\alpha>0$, where $t_N=\ceil{(2\gamma)^{-1}N\log N}$,
\[
\Eb\frac{\xi_v(T_{N,\alpha})}{N}\ge \frac{\rho}{Q}+(1+o(1))\frac{c_0e^\alpha}{\sqrt{N}}
\ge \frac{\rho}{Q}+\frac{c_0 e^\alpha}{2\sqrt{N}}.
\]
This shows that by the Chebyshev inequality and by \eqref{Eq:load_var},
\begin{align*}
&\Pb_{\xi(0)}\(\frac{\xi_v(T_{N,\alpha})}{N}\le \frac{\rho}{Q}+\frac{c_0 e^\alpha}{4\sqrt{N}}\)\\
&\le \Pb_{\xi(0)}\(\frac{\xi_v(T_{N,\alpha})}{N}-\Eb\frac{\xi_v(T_{N,\alpha})}{N} \le -\frac{c_0 e^\alpha}{4\sqrt{N}}\)
\le \frac{16N\Var(\xi_v(T_{N,\alpha})/N)}{c_0^2e^{2\alpha}}= O(e^{-2\alpha}).
\end{align*}

Let us analyze the stationary distribution $\pi_n$.
By \eqref{Eq:load_eq}, \eqref{Eq:load_var} and by taking $t\to \infty$,
the expectation and the variance of $\xi_v/N$ under $\pi_n$ satisfy that for all $N>0$,
\[
\left|\Eb_{\pi_n}\frac{\xi_v}{N}-\frac{\rho}{N}\right|\le \frac{C}{N}
\quad \text{and} \quad
\Var_{\pi_n}\(\frac{\xi_v}{N}\) \le \frac{C}{\gamma N}.
\]
The Chebyshev inequality yields for a constant $C>0$ and for all $R>0$,
\begin{align*}
\pi_n\(\left|\frac{\xi_v}{N}-\frac{\rho}{Q}\right| \ge \frac{R}{\sqrt{N}}\)
&\le \pi_n\(\left|\frac{\xi_v}{N}-\Eb_{\pi_n}\frac{\xi_v}{N}\right|\ge \frac{R}{\sqrt{N}}-\frac{C}{N}\)\\
&\le \frac{\Var_{\pi_n}(\xi_v/N)}{(R/\sqrt{N}-C/N)^2}=O\(\frac{1}{R^2}\).
\end{align*}

Summarizing the above estimates,
we obtain
\[
\Pb_{\xi(0)}\(\frac{\xi_v(T_{N,\alpha})}{N} \ge \frac{\rho}{Q}+\frac{c_0 e^\alpha}{4\sqrt{N}}\)
-\pi_n\(\frac{\xi_v}{N} \ge \frac{\rho}{Q}+\frac{c_0 e^\alpha}{4\sqrt{N}}\)=1-O(e^{-2\alpha}).
\]
Thus $D_\TV(T_{N, \alpha})\ge 1-O(e^{-2\alpha})$,
and letting $N \to \infty$ and then $\alpha \to \infty$ yields the claim.
\qed

\proof[Proof of Theorem \ref{Thm:load}]
The claim follows from Theorems \ref{Thm:load_upper} and \ref{Thm:load_lower}.
\qed

\section{The congestion dynamics}\label{Sec:congestion}

For each $\sigma \in \Sc_n$ and for each $u \in S$,
recall that
\[
\xi(\sigma)=\(\xi_u(\sigma)\)_{u\in S} \quad \text{where $\xi_u(\sigma)=\#\{i \in [n] \ : \ \sigma_i=u\}$}.
\]
For all $\sigma_0, \sigma \in \Sc_n$,
let us define
\[
\xi^{\sigma_0}(\sigma):=\(\xi_{u, v}^{\sigma_0}(\sigma)\)_{u, v \in S}
\quad \text{where $\xi_{u, v}^{\sigma_0}(\sigma):=\#\left\{i \in [n] \ : \ \sigma_{0i}=u, \sigma_i=v\right\}$}.
\]
For the chain $\{\sigma(t)\}_{t\in \Z_+}$ on $\Sc_n$ with $\sigma(0)=\sigma_0$,
let
\[
\xi^{\sigma_0}(t):=\(\xi_{u, v}^{\sigma_0}(t)\)_{u, v \in S} \quad \text{where $\xi_{u, v}^{\sigma_0}(t):=\xi_{u, v}^{\sigma_0}(\sigma(t))$},
\]
and call $\{\xi^{\sigma_0}(t)\}_{t\in \Z_+}$ the {\bf load matrix chain} on $\{0, 1, \dots, \kappa\}^{S\times S}$.
Each $\xi_{u, v}^{\sigma_0}(t)$ counts the number of players who adapt $u$ at time $0$ and adapt $v$ at time $t$.
Furthermore,
for each fixed $\sigma_0 \in \Sc_n$,
let $\wbar \mu_n^{\sigma_0}$ be the law of $\xi^{\sigma_0}(\sigma)$ where $\sigma$ is distributed according to $\mu_n$.

\begin{lemma}\label{Lem:TV_stationary}
For all positive integers $n$, all $\sigma_0 \in \Sc_n$ and $t \in \Z_+$,
it holds that
\[
\|\Pb_{\sigma_0}\(\sigma(t) \in \cdot\)-\mu_n\|_\TV=\|\Pb_{\sigma_0}\(\xi^{\sigma_0}(t)\in \cdot\)-\wbar \mu_n^{\sigma_0}\|_\TV.
\]
\end{lemma}

\proof
For each $\sigma_0\in \Sc_n$,
let
\[
\Mcc(\sigma_0):=\left\{\(\xi_{u, v}\)_{u, v \in S} \in \{0, 1, \dots, \kappa\}^{S\times S} \ : \ \(\sum_{v\in S}\xi_{u, v}\)_{u \in S}=\xi(\sigma_0)\right\}.
\]
Furthermore for each $\mb{\xi} \in \Mcc(\sigma_0)$,
let
\[
\Sc(\sigma_0, \mb{\xi}):=\Big\{\sigma \in \Sc_n \ : \ \xi^{\sigma_0}(\sigma)=\mb{\xi}\Big\}.
\]
Note that for each $\sigma_0 \in \Sc_n$,
we have $\Sc_n=\bigsqcup_{\mb{\xi}\in \Mcc(\sigma_0)}\Sc(\sigma_0, \mb{\xi})$.
The proof uses the following observation:
For each $\sigma_0 \in \Sc_n$ and each $t \in \Z_+$,
the distribution of $\sigma(t)$ with $\sigma(0)=\sigma_0$ is invariant under permutations preserving $\sigma_0$ on $[n]$.
This implies that for each $\sigma \in \Sc(\sigma_0, \mb{\xi})$,
\[
\Pb_{\sigma_0}\(\sigma(t)=\sigma \mid \sigma(t)\in \Sc(\sigma_0, \mb{\xi})\)=\frac{1}{\#\Sc(\sigma_0, \mb{\xi})}.
\]
Moreover, since the stationary distribution $\mu_n$ is invariant under all permutations on $[n]$,
in particular,
all those preserving $\sigma_0$,
for each $\sigma \in \Sc(\sigma_0, \mb{\xi})$,
\[
\mu_n\(\sigma \mid \Sc(\sigma_0, \mb{\xi})\)=\frac{1}{\#\Sc(\sigma_0, \mb{\xi})}.
\]
Therefore it follows that for each $\sigma_0 \in \Sc_n$ and each $t \in \Z_+$,
\begin{align*}
&\|\Pb_{\sigma_0}(\sigma(t) \in \cdot)-\mu_n\|_\TV
=\frac{1}{2}\sum_{\mb{\xi}\in \Mcc(\sigma_0)}\sum_{\sigma \in \Sc(\sigma_0, \mb{\xi})}|\Pb_{\sigma_0}(\sigma(t)=\sigma)-\mu_n(\sigma)|\\
&=\frac{1}{2}\sum_{\mb{\xi} \in \Mcc(\sigma_0)}\sum_{\sigma\in \Sc(\xi_0, \mb{\xi})}\frac{1}{\#\Sc(\sigma_0, \mb{\xi})}|\Pb_{\sigma_0}(\xi^{\sigma_0}(t)=\mb{\xi})-\mu_n\(\{\sigma \ : \ \xi^{\sigma_0}(\sigma)=\mb{\xi}\}\)|.
\end{align*}
This shows that
\[
\|\Pb_{\sigma_0}(\sigma(t)\in \cdot)-\mu_n\|_\TV=\frac{1}{2}\sum_{\mb{\xi}\in \Mcc(\sigma_0)}|\Pb_{\sigma_0}(\xi^{\sigma_0}(t)=\mb{\xi})-\wbar \mu_n^{\sigma_0}(\mb{\xi})|,
\]
concluding the claim.
\qed

\medskip

For $\delta>0$,
let
\[
\Sc(\delta):=\left\{\sigma \in \Sc_n \ : \ \left\|\frac{\xi(\sigma)}{N}-\(\frac{\rho}{Q}\)_{u\in S}\right\| \le \delta \right\}.
\]
In the notation in Section \ref{Sec:load},
we have $\Sc(\delta)=\left\{\sigma \in \Sc_n \ : \ \xi(\sigma) \in \Xi(\delta)\right\}$.

\begin{lemma}\label{Lem:burning-in}
Let $T_{\alpha, N}:=\ceil{\alpha N}$ for reals $\alpha$ and integers $N$.
For each $\delta>0$ there exist constants $c_\delta, \alpha_\delta, N_\delta>0$
such that for all $\alpha \ge \alpha_\delta$, 
$N\ge N_\delta$,
and $\sigma_\ast \in \Sc_n$,
\[
\Pb_{\sigma_\ast}\(\sigma(T_{\alpha, N})\notin \Sc(\delta)\)\le e^{-c_\delta N}.
\]
Furthermore, if $\sigma_\ast \in \Sc(\delta)$,
then for all $N\ge N_\delta$,
\[
\Pb_{\sigma_\ast}\Big(\bigcup_{0\le t\le N^2}\big\{\sigma(t) \notin \Sc(2\delta)\big\}\Big)\le e^{-c_\delta N}.
\]
\end{lemma}

\proof
Applying to Lemma \ref{Lem:general} by letting $x_t=\|\Delta(t)\|^2$
with $\lambda=2\gamma$ in \eqref{Eq:load} implies that there exist constants $c, R_0>0$ such that for all $R \ge R_0$ and for each $t \in \Z_+$,
\[
\Pb_{\sigma_\ast}\(\|\Delta(t)\|^2\ge \(1-\frac{2\gamma}{N}\)^t \|\Delta(0)\|^2+\frac{R}{\sqrt{N}}\) \le e^{-c R^2}.
\]
Hence since $\|\Delta(0)\|\le \sqrt{Q}$ and $(1-2\gamma/N)^{T_{\alpha, N}}\le e^{-2\gamma \alpha}$,
for $\delta>0$ and for $R=\delta^2 \sqrt{N}/2$,
\[
\Pb_{\sigma_\ast}\(\|\Delta(T_{\alpha, N})\|^2 \ge Q e^{-2\gamma \alpha}+\frac{\delta^2}{2}\) \le e^{-c\delta^4N/4}.
\]
For each $\delta>0$ we have $Q e^{-2\gamma \alpha} \le \delta^2/2$ for all large enough $\alpha$.
This concludes the first clam.
The second claim follows from Lemma \ref{Lem:load_good}.
\qed

\medskip

For each $\sigma_0 \in \Sc_n$ and each $t \in \Z_+$,
we find that for $u, v \in S$,
\[
\xi_{u,v}^{\sigma_0}(t+1)-\xi_{u,v}^{\sigma_0}(t)=\frac{\xi_u(\sigma_0)-\xi_{u,v}^{\sigma_0}(t)}{n}\frac{\kappa-\xi_v(t)}{1+N-n}-\frac{\xi_{u,v}^{\sigma_0}(t)}{n}\frac{N-n-(\kappa-\xi_v(t))}{1+N-n}+M_{u,v,t+1},
\]
where $|M_{u,v,t+1}|\le 2$ and $\Eb[M_{u,v,t+1}\mid \Fc_t]=0$ almost surely.
Rearranging terms,
we obtain
\[
\xi_{u,v}^{\sigma_0}(t+1)-\xi_{u,v}^{\sigma_0}(t)=\frac{\xi_u(\sigma_0)}{n}\frac{\kappa-\xi_v(t)}{1+N-n}-\frac{\xi_{u, v}^{\sigma_0}(t)}{n}\frac{N-n}{1+N-n}+M_{u,v,t+1}.
\]
It is convenient to write
\begin{equation}\label{Eq:load_matrix}
\xi_{u,v}^{\sigma_0}(t+1)-\xi_{u,v}^{\sigma_0}(t)=\frac{\xi_u(\sigma_0)}{n}\frac{\kappa-\xi_v(t)}{(1-\rho)N}-\frac{\xi_{u,v}^{\sigma_0}(t)}{n}+\wh M_{u,v,t+1},
\end{equation}
where $|\wh M_{u,v,t+1}|\le C$ and $|\Eb[\wh M_{u,v,t+1}\mid \Fc_t]|\le C/n$ almost surely for a constant $C$.
If $\xi_u(\sigma_0)>0$ for every $u \in S$,
then let
\[
\wbar \xi^{\sigma_0}(t):=\(\wbar \xi^{\sigma_0}_{u, v}(t)\)_{u,v \in S} 
\quad \text{where $\wbar \xi^{\sigma_0}_{u,v}(t):=\frac{\xi^{\sigma_0}_{u,v}(t)}{\xi_u(\sigma_0)}$}.
\]
In the case when these are defined,
we have
\begin{equation}\label{Eq:load_matrix_eq}
\wbar \xi^{\sigma_0}_{u,v}(t+1)-\wbar \xi^{\sigma_0}_{u,v}(t)=\frac{\kappa-\xi_v(t)}{(1-\rho)nN}-\frac{\wbar \xi^{\sigma_0}_{u,v}(t)}{n}+\frac{\wh M_{u,v,t+1}}{\xi_u(\sigma_0)}.
\end{equation}

For an arbitrary $\sigma_0 \in \Sc_n$ such that $\xi_u(\sigma_0)>0$ for all $u \in S$,
let us define for $\sigma \in \Sc_n$,
\[
\mb{\Delta}_u^{\sigma_0}(\sigma):=\(\Delta_{u,v}^{\sigma_0}(\sigma)\)_{v \in S}
\quad \text{for $u\in S$ where $\Delta_{u,v}^{\sigma_0}(\sigma):=\wbar \xi^{\sigma_0}_{u,v}(\sigma)-\frac{1}{Q}$}.
\]
Note that if $0<\delta<\rho/(2Q)$ and $\sigma_0 \in \Sc(\delta)$,
then $\wbar \xi^{\sigma_0}_{u,v}(\sigma)$ is well-defined for all $\sigma \in \Sc_n$.
Simplifying notations, whenever $\wbar \xi^{\sigma_0}_{u,v}(\sigma(t))$ is defined,
we write for $\{\sigma(t)\}_{t\in \Z_+}$,
\[
\mb{\Delta}_u^{\sigma_0}(t):=\(\Delta_{u,v}^{\sigma_0}(t)\)_{v \in S}
\quad \text{for $u \in S$ where $\Delta^{\sigma_0}_{u,v}(t):=\Delta^{\sigma_0}_{u,v}(\sigma(t))$}.
\]
For all $u, v \in S$ and each $t \in \Z_+$,
it holds that by \eqref{Eq:load_matrix_eq},
\[
\Delta^{\sigma_0}_{u,v}(t+1)=\(1-\frac{1}{n}\)\Delta^{\sigma_0}_{u,v}(t)+\frac{1}{n}\(\frac{\kappa-\xi_v(t)}{(1-\rho)N}-\frac{1}{Q}\)+\frac{\wh M_{u,v,t+1}}{\xi_u(\sigma_0)}.
\]
Let us write $\mb{\kappa}:=(\kappa, \dots, \kappa)$ and $1/\mb{Q}:=(1/Q, \dots, 1/Q)$ in $\R^Q$.
For each $u \in S$,
we have the following $\R^Q$-valued process:
\begin{equation*}
\mb{\Delta}^{\sigma_0}_u(t+1)=\(1-\frac{1}{n}\)\mb{\Delta}^{\sigma_0}_u(t)+\frac{1}{n}\(\frac{1}{(1-\rho)N}(\mb{\kappa}-\xi(t))-\frac{1}{\mb{Q}}\)+\mb{M}_{u,t+1},
\end{equation*}
where
$\mb{M}_{u, t+1}:=\(\wh M_{u,v,t+1}/\xi_u(\sigma_0)\)_{v \in S}$.
Note that 
\begin{equation*}
\frac{1}{(1-\rho)N}\(\mb{\kappa}-\xi(t)\)-\frac{1}{\mb{Q}}=-\frac{1}{1-\rho}\(\frac{\xi(t)}{N}-\(\frac{\rho}{Q}\)_{v\in S}\)=-\frac{1}{1-\rho}\Delta(t).
\end{equation*}
Therefore we have
\begin{equation}\label{Eq:Delta_matrix}
\mb{\Delta}^{\sigma_0}_u(t+1)=\(1-\frac{1}{n}\)\mb{\Delta}^{\sigma_0}_u(t)-\frac{1}{n(1-\rho)}\Delta(t)+\mb{M}_{u,t+1},
\end{equation}

For $0<\delta<\rho/(2Q)$ and $\sigma_0 \in \Sc(\delta)$,
let
\[
\Sc(\sigma_0, \delta):=\left\{\sigma \in \Sc_n \ : \ \max_{u \in S}\|\mb{\Delta}^{\sigma_0}_u(\sigma)\| \le \delta \right\}.
\]

\begin{lemma}\label{Lem:load_matrix_good}
There exists a constant $c_\star>0$ depending only on $Q$ and $\rho$ satisfying the following:
For all $0<\delta<\rho/(2Q)$
there exist $c_\delta, n_\delta>0$ such that if $\sigma_0 \in \Sc(\delta)$ and $\sigma_\ast \in \Sc(\sigma_0, \delta)$,
then for all $n \ge n_\delta$,
\[
\Pb_{\sigma_\ast}\Bigg(\bigcup_{0\le t \le n^2}\Big\{\sigma(t) \notin \Sc\(\sigma_0, c_\star\sqrt{\delta}\)\Big\}\Bigg)\le e^{-c_\delta n}.
\]
\end{lemma}

\proof
Note that if $0<\delta<\rho/(2Q)$ and $\sigma_0 \in \Sc(\delta)$,
then for each $u \in S$,
\begin{equation}\label{Eq:sigma_0}
\xi_u(\sigma_0) \ge \(\frac{\rho}{Q}-\delta\)N>\frac{\rho}{2Q}N\ge \frac{n}{2Q}.
\end{equation}
Further if $\sigma_\ast \in \Sc(\sigma_0, \delta)$,
then $\sigma_\ast \in \Sc(\delta)$.
Indeed, if $\sigma_\ast \in \Sc(\sigma_0, \delta)$, 
then
\[
\left\|(\xi_{u, v}^{\sigma_0}(\sigma_\ast))_{v \in S} -(\xi_u(\sigma_0)/Q)_{v \in S}\right\|\le \delta \xi_u(\sigma_0)
\]
for each $u \in S$, and thus
summing over $u \in S$ and dividing by $N$ shows that $\sigma_\ast\in \Sc(\delta)$.

By \eqref{Eq:Delta_matrix},
squaring the norms of both sides yields
\begin{equation}\label{Eq:Lem:load_matrix_good:0}
\|\mb{\Delta}^{\sigma_0}_u(t+1)\|^2=\left\|\(1-\frac{1}{n}\)\mb{\Delta}^{\sigma_0}_u(t)-\frac{1}{n(1-\rho)}\Delta(t)\right\|^2+\wt M_{u, t+1},
\end{equation}
where $|\wt M_{u, t+1}|\le CQ^2/n$ and $|\Eb[\wt M_{u, t+1}\mid \Fc_t]|\le CQ^2/n^2$ almost surely by \eqref{Eq:sigma_0}.

Let us define the event for $t \in \Z_+$ and $\delta>0$,
\[
\Ac_{t, \delta}:=\bigcap_{0\le s \le t}\Big\{\sigma(s) \in \Sc(2\delta)\Big\}.
\]
By Lemma \ref{Lem:load_good},
since $\sigma_\ast \in \Sc(\delta)$,
we have
\begin{equation}\label{Eq:Lem:load_matrix_good:1}
\Pb_{\sigma_\ast}\Bigg(\bigcup_{0\le t \le n^2}\Ac_{t, \delta}^{\sf c}\Bigg) \le e^{-c_\delta n}.
\end{equation}
Since $\Ac_{t+1, \delta} \subset \Ac_{t, \delta}$ and $\Ac_{t, \delta} \in \Fc_t$,
we have by \eqref{Eq:Lem:load_matrix_good:0} and by letting $M'_{u, t+1}:=\wt M_{u, t+1}\1_{\Ac_{t, \delta}}$,
\[
\|\mb{\Delta}^{\sigma_0}_u(t+1)\|^2\1_{\Ac_{t+1,\delta}}
\le \left\|\(1-\frac{1}{n}\)\mb{\Delta}^{\sigma_0}_u(t)-\frac{1}{n(1-\rho)}\Delta(t)\right\|^2\1_{\Ac_{t, \delta}}+M_{u, t+1}'.
\]
Furthermore, we obtain by using $\|\Delta(t)\|\le 2\delta$ on $\Ac_{t, \delta}$ and $\|\mb{\Delta}^{\sigma_0}_u(t)\| \le \sqrt{Q}$, 
\begin{align*}
&\|\mb{\Delta}^{\sigma_0}_u(t+1)\|^2\1_{\Ac_{t+1,\delta}}\\
&\qquad \qquad\le \(1-\frac{1}{n}\)^2\|\mb{\Delta}^{\sigma_0}_u(t)\|^2\1_{\Ac_{t, \delta}}+\frac{1}{n^2(1-\rho)^2}\|\Delta(t)\|^2+\frac{2\delta \sqrt{Q}}{(1-\rho)n}+M_{u, t+1}'.
\end{align*}
Rearranging the terms and using $\|\mb{\Delta}_u^{\sigma_0}(t)\|^2\1_{\Ac_{t, \delta}}\le Q$ and $\|\Delta(t)\|^2\le Q$,
we have
\begin{align*}
&\|\mb{\Delta}_u^{\sigma_0}(t+1)\|^2\1_{\Ac_{t+1, \delta}}-\frac{\delta \sqrt{Q}}{1-\rho}\\
&\qquad \qquad\le \(1-\frac{2}{n}\)\(\|\mb{\Delta}_u^{\sigma_0}(t)\|^2\1_{\Ac_{t,\delta}}-\frac{\delta \sqrt{Q}}{1-\rho}\)+\frac{Q}{n^2(1-\rho)^2}+\frac{Q}{n^2}+M'_{u, t+1}.
\end{align*}
Applying to Lemma \ref{Lem:general} by setting $x_t=\|\mb{\Delta}_u^{\sigma_0}(t)\|^2\1_{\Ac_{t,\delta}}-\delta \sqrt{Q}/(1-\rho)$ and $\lambda=2$,
we have $c>0$ and $R_0$ such that for all $R \ge R_0$, $t \in \Z_+$ and $n>0$, 
\[
\Pb_{\sigma_\ast}\(\|\mb{\Delta}_u^{\sigma_0}(t)\|^2\1_{\Ac_{t, \delta}}\ge \frac{\delta \sqrt{Q}}{1-\rho}+\frac{R}{\sqrt{n}}+\(1-\frac{2}{n}\)^t x_0\)\le e^{-cR^2}.
\]
If $\|\mb{\Delta}_u^{\sigma_0}(0)\|^2 \le \delta \sqrt{Q}/(1-\rho)$,
then $x_0 \le 0$ and thus by setting $R=\delta^2 \sqrt{n}$,
\[
\Pb_{\sigma_\ast}\(\|\mb{\Delta}_u^{\sigma_0}(t)\|^2\1_{\Ac_{t, \delta}} \ge \frac{\delta \sqrt{Q}}{1-\rho}+\delta^2\)\le e^{-c\delta^4 n}.
\]
Let us define the constant
\[
c_{Q, \rho}:=\frac{\sqrt{Q}}{1-\rho}+\frac{\rho}{2Q}.
\]
Since $0< \delta <\rho/(2Q)$,
we have
$c_{Q, \rho}\delta>\delta \sqrt{Q}(1-\rho)^{-1}+\delta^2$.
Thus, if $\|\mb{\Delta}_u^{\sigma_0}(0)\|^2 \le \delta^2$, 
then by the union bound over $0 \le t \le n^2$, for all $n>0$,
\[
\Pb_{\sigma_\ast}\(\max_{0 \le t\le n^2}\|\mb{\Delta}_u^{\sigma_0}(t)\|^2\1_{\Ac_{t, \delta}} \ge c_{Q, \rho}\delta\)
\le n^2 e^{-c\delta^4 n} \le e^{-c_\delta n}.
\]
By \eqref{Eq:Lem:load_matrix_good:1},
if $\sigma_0 \in\Sc(\delta)$ and $\sigma_\ast \in \Sc(\sigma_0, \delta)$,
then for all $n>0$, 
\begin{align*}
&\Pb_{\sigma_\ast}\(\max_{0\le t\le n^2}\|\mb{\Delta}_u^{\sigma_0}(t)\|^2 \ge c_{Q, \rho}\delta\)\\
&\le\Pb_{\sigma_\ast}\(\bigcup_{0 \le t \le n^2}\left\{\|\mb{\Delta}_u^{\sigma_0}(t)\|^2\1_{\Ac_{t, \delta}}\ge c_{Q,\rho}\delta\right\}\cap \Ac_{t,\delta}\)+\Pb_{\sigma_\ast}\(\bigcup_{0\le t\le n^2}\Ac_{t, \delta}^{\sf c}\)\le 2 e^{-c_\delta n}.
\end{align*}
Further taking the union bound over $u \in S$ shows the following:
For $c_\star:=\sqrt{c_{Q, \rho}}$,
for all $0<\delta<\rho/(2Q)$,
if $\sigma_0 \in \Sc(\delta)$ and $\sigma_\ast \in \Sc(\sigma_0, \delta)$,
\[
\Pb_{\sigma_\ast}\(\max_{u \in S}\max_{0\le t\le n^2}\|\mb{\Delta}_u^{\sigma_0}(t)\| \ge c_\star\sqrt{\delta}\)\le 2Q e^{-c_\delta n}.
\]
This yields the claim.
\qed

\begin{lemma}\label{Lem:load_matrix_square_root}
Let $t_n:=\ceil{(1/2)n\log n}$.
There exists a constant $C$ such that the following holds for all large enough $n$:
For all $0<\delta<\rho/(2Q)$, 
let $\sigma_0 \in \Sc(\delta)$.
For all $\sigma_\ast \in \Sc_n$ and
$R>1$,
\[
\Pb_{\sigma_\ast}\(\sigma(t_n) \notin \Sc\(\sigma_0, \frac{R}{\sqrt{n}}\)\) \le \frac{CQ^2}{R}.
\]
\end{lemma}

\proof
For all $0<\delta<\rho/(2Q)$,
if $\sigma_0\in \Sc(\delta)$,
then $\xi_u(\sigma_0) \ge N\rho/(2Q)>0$ for all $u \in S$ and $\wbar \xi^{\sigma_0}(t)$ is well-defined for all $t \in \Z_+$.
Let us define the following function on $\R^Q$ for each positive integer $n$,
\[
\f_n(x):=\|x\|+\frac{1}{\sqrt{n}}e^{-\sqrt{n}\|x\|}-\frac{1}{\sqrt{n}}.
\]
It holds that for all $x, h \in \R^Q$,
\[
\f_n(x+h) \le \f_n(x)+\abr{h, \nabla \f_n(x)}+\frac{\sqrt{n}}{2}\|h\|^2 \quad \text{and} \quad \|\nabla \f_n(x)\| \le 1,
\]
where $\abr{\cdot,\cdot}$ is the standard inner product in $\R^Q$ \cite[Lemma 3.6]{PeresTanakaZhai}.
By \eqref{Eq:Delta_matrix},
\[
\f_n\(\mb{\Delta}^{\sigma_0}_u(t+1)\)\le \f_n\(\(1-\frac{1}{n}\)\mb{\Delta}^{\sigma_0}_u(t)-\frac{1}{n(1-\rho)}\Delta(t)\)+\wt M_{u, t+1},
\]
where $|\wt M_{u, t+1}|\le CQ/n$ and $|\Eb[\wt M_{u, t+1}\mid \Fc_t]|\le CQ/n^{3/2}$ almost surely.
Hence for each $u \in S$,
\begin{equation}\label{Eq:Lem:Delta_phi}
\f_n(\mb{\Delta}_u^{\sigma_0}(t+1))\le \(1-\frac{1}{n}\)\f_n\(\mb{\Delta}_u^{\sigma_0}(t)\)
+\frac{1}{n(1-\rho)}\left\|\Delta(t)\right\|+\wt M_{u, t+1}.
\end{equation}
This follows since $\f_n$ is convex and $\f_n(x)\le \|x\|$ for all $x \in \R^Q$.
By \eqref{Eq:load_eq},
almost surely
for each $t \in \Z_+$,
\[
\Eb\left[\left\|\Delta(t+1)\right\|^2 \mid\Fc_t \right]
\le \(1-\frac{\gamma}{N}\)^2\left\|\Delta(t)\right\|^2+\frac{CQ}{N^2}.
\]
This implies that by induction in $t \in \Z_+$ and by $\|\Delta(0)\|^2\le Q$,
\[
\Eb\|\Delta(t)\|^2 \le \(1-\frac{\gamma}{N}\)^{2t}\Eb\|\Delta(0)\|^2+\frac{C'Q}{N}\le \(1-\frac{\gamma}{N}\)^{2t}Q+\frac{C'Q}{N}.
\]
Hence by the Cauchy-Schwarz inequality and the inequality $\sqrt{t_1+t_2}\le \sqrt{t_1}+\sqrt{t_2}$ for $t_i\ge 0$, $i=1, 2$,
\[
\Eb\|\Delta(t)\|\le \(1-\frac{\gamma}{N}\)^t\sqrt{Q}+\sqrt{\frac{C'Q}{N}}.
\]
Taking the expectations in \eqref{Eq:Lem:Delta_phi},
we obtain for each $u \in S$ and each $t \in \Z_+$,
\[
\Eb\f_n(\mb{\Delta}^{\sigma_0}_u(t+1))\le \(1-\frac{1}{n}\)\Eb \f_n(\mb{\Delta}^{\sigma_0}_u(t))+\(1-\frac{\gamma}{N}\)^t\frac{\sqrt{Q}}{n(1-\rho)}+\frac{C''Q}{n^{3/2}}.
\]
By induction in $t \in \Z_+$,
it follows that
\begin{align*}
&\Eb \f_n(\mb{\Delta}_u^{\sigma_0}(t))\\
&\le \(1-\frac{1}{n}\)^t\Eb \f_n(\mb{\Delta}^{\sigma_0}_u(0))+
\frac{\sqrt{Q}}{n(1-\rho)}\sum_{s=0}^{t-1}\(1-\frac{1}{n}\)^s\(1-\frac{\gamma}{N}\)^{t-1-s}+\frac{CQ}{\sqrt{n}}.
\end{align*}
Note that for all large enough $n$,
we have $n=\rho N+O(1)$ and
\[
1-\frac{1}{n} >1-\frac{1}{\rho(1-\rho)N}=1-\frac{\gamma}{N}.
\]
Thus we have
\begin{align*}
&\sum_{s=0}^{t-1}\(1-\frac{1}{n}\)^s\(1-\frac{\gamma}{N}\)^{t-1-s}
=\sum_{s=0}^{t-1}\(1-\frac{1}{n}\)^{t-1}\(\frac{1-\gamma/N}{1-1/n}\)^{t-1-s}\\
&\le \(1-\frac{1}{n}\)^{t-1}\frac{1}{1-(1-\gamma/N)/(1-1/n)}
=\(1-\frac{1}{n}\)^t\frac{N}{\gamma-N/n}.
\end{align*}
Summarizing the above estimates,
we obtain for some constant $C>0$,
for all large enough $n$ and for all $t \in \Z_+$,
\[
\Eb\f_n(\mb{\Delta}^{\sigma_0}_u(t))\le \(1-\frac{1}{n}\)^t\Eb\f_n(\mb{\Delta}^{\sigma_0}_u(0))+C\sqrt{Q}\(1-\frac{1}{n}\)^t+\frac{CQ}{\sqrt{n}}.
\]
Therefore for $t_n=\ceil{(1/2)n\log n}$,
\[
\Eb \f_n(\mb{\Delta}^{\sigma_0}_u(t_n))\le \frac{1}{\sqrt{n}}\Eb\f_n(\mb{\Delta}^{\sigma_0}_u(0))+\frac{C\sqrt{Q}}{\sqrt{n}}+\frac{CQ}{\sqrt{n}}\le \frac{C'Q}{\sqrt{n}},
\]
where we have used $\Eb \f_n(\mb{\Delta}^{\sigma_0}_u(0))\le \Eb\|\mb{\Delta}^{\sigma_0}_u(0)\|\le \sqrt{Q}$.
It follows that for all $R>1$,
\begin{align*}
&\Pb_{\sigma_\ast}\(\max_{u \in S}\|\mb{\Delta}_u^{\sigma_0}(t_n)\|>\frac{R}{\sqrt{n}}\)
\le \Pb_{\sigma_\ast}\(\max_{u\in S}\f_n(\mb{\Delta}^{\sigma_0}_u(t_n))>\frac{R-1}{\sqrt{n}}\)\\
&\le \sum_{u \in S}\Pb_{\sigma_\ast}\(\f_n(\mb{\Delta}^{\sigma_0}_u(t_n))>\frac{R-1}{\sqrt{n}}\)
\le \sum_{u \in S}\frac{\sqrt{n}}{R-1}\Eb \f_n(\mb{\Delta}^{\sigma_0}_u(t_n))\le \frac{CQ^2}{R}.
\end{align*}
In the above the Markov inequality has been used in the third inequality.
This concludes the claim.
\qed

\medskip

For each $\sigma_0 \in \Sc_n$,
let us define for $\sigma, \wt \sigma \in \Sc_n$,
\[
\mb{D}^{\sigma_0}(\sigma, \wt \sigma):=\sum_{u\in S}D_u^{\sigma_0}(\sigma, \wt \sigma) 
\quad \text{where $D_u^{\sigma_0}(\sigma, \wt \sigma):=\frac{1}{2}\sum_{v \in S}|\xi_{u,v}^{\sigma_0}(\sigma)-\xi_{u,v}^{\sigma_0}(\wt \sigma)|$}. 
\]
Let us also define
\[
D(\sigma, \wt \sigma):=\frac{1}{2}\sum_{v \in S}|\xi_v(\sigma)-\xi_v(\wt \sigma)|.
\]
We have $D(\sigma, \wt \sigma)=D(\xi(\sigma), \xi(\wt \sigma))$ in the notation of Section \ref{Sec:load}.
For chains $\{\sigma(t)\}_{t\in \Z_+}$ and $\{\wt \sigma(t)\}_{t\in \Z_+}$,
let $\mb{D}_t^{\sigma_0}:=\mb{D}^{\sigma_0}(\sigma(t), \wt \sigma(t))$
and $D_t=D(\sigma(t), \wt \sigma(t))$ for brevity.

\begin{lemma}\label{Lem:coupling_matrix}
Let $\rho \in (0, 1/2]$.
Fix $\delta>0$ such that
\[
0<\max\Big\{\delta, c_\star \sqrt{\delta}\Big\} \le \frac{\rho}{8Q},
\]
where $c_\star$ is the constant in Lemma \ref{Lem:load_matrix_good}.
There exists a coupling $\Pb$ between the chains $\{\sigma(t)\}_{t\in \Z_+}$ and $\{\wt \sigma(t)\}_{t\in \Z_+}$ satisfying the following:
For each $\sigma_0\in \Sc(\delta)$ and $t \in \Z_+$,
if $\sigma(t), \wt \sigma(t) \in \Sc\(2\delta\) \cap\Sc\big(\sigma_0, c_\star\sqrt{\delta}\big)$ and $\mb{D}_t^{\sigma_0}>0$,
then almost surely for each $r=1$, $-1$,
\[
\Pb\(\mb{D}^{\sigma_0}_{t+1}-\mb{D}^{\sigma_0}_t = r \mid \sigma(t), \wt \sigma(t)\)\ge \frac{1}{4Q^3}
\quad \text{and} \quad
\Eb\left[\mb{D}^{\sigma_0}_{t+1}-\mb{D}^{\sigma_0}_t\mid \sigma(t), \wt \sigma(t)\right] \le 0.
\]
\end{lemma}

\proof
Let us simply write $\sigma$ for $\sigma(t)$, and $\wt \sigma$ for $\wt \sigma(t)$.
If $\sigma, \wt \sigma \in \Sc\(2\delta\)\cap \Sc\big(\sigma_0, c_\star\sqrt{\delta}\big)$ and $\sigma_0\in \Sc(\delta)$,
then the following holds:
For all $u, v \in S$,
\[
\xi_{u,v}^{\sigma_0}(\sigma) \ge \(\frac{1}{Q}-c_\star\sqrt{\delta}\)\xi_u(\sigma_0) \ge \(\frac{1}{Q}-c_\star\sqrt{\delta}\)\(\frac{\rho}{Q}-\delta\)N\ge \frac{\rho N}{2Q^2} \ge \frac{n}{2Q^2}.
\]
Similarly, it holds that for all $v \in S$,
\[
\kappa-\xi_v(\sigma) \ge \kappa-\(\frac{\rho}{Q}+2\delta\)N
=\(\frac{1-\rho}{Q}-2\delta\) N
\ge \frac{(1-\rho)N+1}{2Q}\ge\frac{N-n}{2Q}.
\]
The analogous estimates hold also for $\wt \sigma$.

If $\mb{D}^{\sigma_0}(\sigma, \wt \sigma)>0$,
then there exist $u \in S$ and a pair $(v_0, v_1) \in S\times S$ such that
\[
\xi_{u, v_0}^{\sigma_0}(\sigma)<\xi_{u, v_0}^{\sigma_0}(\wt \sigma) \quad \text{and} \quad \xi_{u, v_1}^{\sigma_0}(\sigma)>\xi_{u, v_1}^{\sigma_0}(\wt \sigma).
\]
Let us fix such $u \in S$ and a pair $(v_0, v_1) \in S \times S$.
For $u, v \in S$,
let
\[
m_{u, v}:=\min\left\{\xi^{\sigma_0}_{u,v}(\sigma), \xi^{\sigma_0}_{u,v}(\wt \sigma)\right\}
\quad \text{and} \quad
M_v:=\max\left\{\xi_v(\sigma), \xi_v(\wt \sigma)\right\}.
\]
Since we assume that $\sigma, \wt \sigma \in \Sc\(2\delta\)\cap \Sc\big(\sigma_0, c_\star\sqrt{\delta}\big)$ and $\sigma_0 \in \Sc(\delta)$,
for all $u, v\in S$,
\begin{equation}\label{Eq:Lem:coupling_matrix:1}
m_{u, v} \ge \ceil{\frac{n}{2Q^2}} \quad \text{and} \quad \kappa-M_v \ge \ceil{\frac{N-n}{2Q}}.
\end{equation}

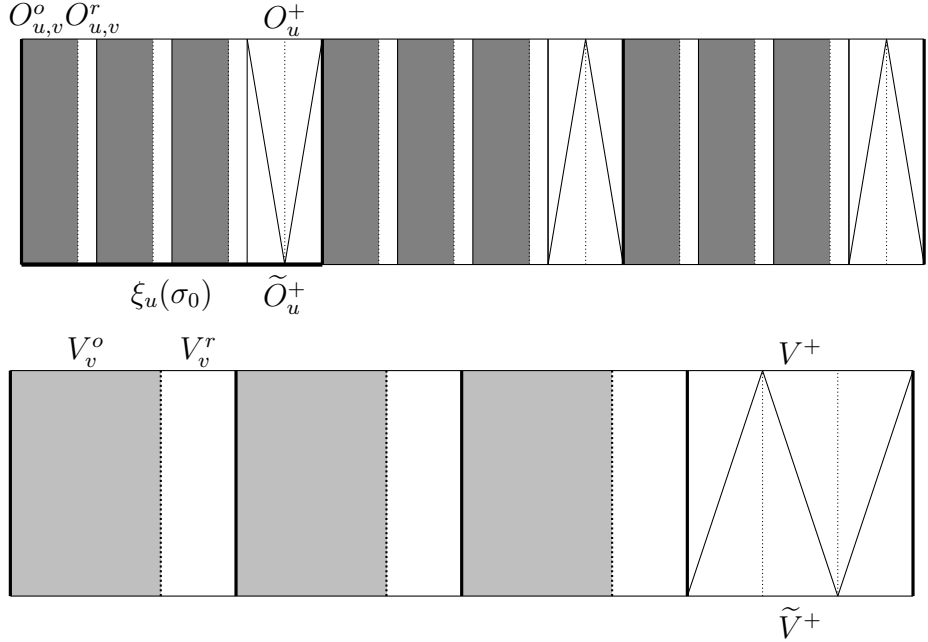
\begin{figure}
\centering

\begin{tikzpicture}[scale=0.5]


\fill [fill=gray] (0,0)--(1.5,0)--(1.5,6)--(0,6)--cycle;
\fill [fill=gray] (2,0)--(3.5,0)--(3.5,6)--(2,6)--cycle;
\fill [fill=gray] (4,0)--(5.5,0)--(5.5,6)--(4,6)--cycle;

\fill [fill=gray] (8,0)--(9.5,0)--(9.5,6)--(8,6)--cycle;
\fill [fill=gray] (10,0)--(11.5,0)--(11.5,6)--(10,6)--cycle;
\fill [fill=gray] (12,0)--(13.5,0)--(13.5,6)--(12,6)--cycle;

\fill [fill=gray] (16,0)--(17.5,0)--(17.5,6)--(16,6)--cycle;
\fill [fill=gray] (18,0)--(19.5,0)--(19.5,6)--(18,6)--cycle;
\fill [fill=gray] (20,0)--(21.5,0)--(21.5,6)--(20,6)--cycle;



\draw (0.4,6.5) node {{$O_{u,v}^o$}};
\draw (1.9,6.5) node {{$O_{u,v}^r$}};
\draw (7,6.5) node {{$O_u^+$}};
\draw (7,-0.8) node {{$\wt O_u^+$}};

\draw[very thick] (0,0)--(0,6);
\draw[densely dotted] (1.5,0)--(1.5,6);
\draw[solid] (2,0)--(2,6);

\draw[densely dotted] (3.5,0)--(3.5,6);
\draw[solid] (4,0)--(4,6);

\draw[densely dotted] (5.5,0)--(5.5,6);
\draw[solid] (6,0)--(6,6);
\draw[solid] (7,0)--(6,6);
\draw[solid] (7,0)--(8,6);
\draw[densely dotted] (7,0)--(7,6);

\draw[very thick] (8,0)--(8,6);
\draw[densely dotted] (9.5,0)--(9.5,6);
\draw[solid] (10,0)--(10,6);

\draw[densely dotted] (11.5,0)--(11.5,6);
\draw[solid] (12,0)--(12,6);

\draw[densely dotted] (13.5,0)--(13.5,6);
\draw[solid] (14,0)--(14,6);

\draw[solid] (14,0)--(14,6);
\draw[solid] (14,0)--(15,6);
\draw[solid] (16,0)--(15,6);
\draw[densely dotted] (15,0)--(15,6);

\draw[very thick] (16,0)--(16,6);
\draw[densely dotted] (17.5,0)--(17.5,6);
\draw[solid] (18,0)--(18,6);

\draw[densely dotted] (19.5,0)--(19.5,6);
\draw[solid] (20,0)--(20,6);

\draw[densely dotted] (21.5,0)--(21.5,6);
\draw[solid] (22,0)--(22,6);

\draw[solid] (22,0)--(22,6);
\draw[solid] (22,0)--(23,6);
\draw[solid] (24,0)--(23,6);
\draw[densely dotted] (23,0)--(23,6);

\draw[very thick] (24,0)--(24,6);

\draw[solid] (0,6)--(24,6);
\draw[solid] (0,0)--(24,0);

\draw[ultra thick] (0,0)--(8,0);
\draw (4,-0.8) node {{$\xi_u(\sigma_0)$}};
\end{tikzpicture}


\begin{tikzpicture}[scale=0.5]


\draw (2,6.5) node {{$V_v^o$}};
\draw (5,6.5) node {{$V_v^r$}};
\draw (21,6.5) node {{$V^+$}};
\draw (21,-0.7) node {{$\wt V^+$}};


\fill [fill=lightgray] (0,0)--(0,6)--(4,6)--(4,0)--cycle;
\fill [fill=lightgray] (6,0)--(6,6)--(10,6)--(10,0)--cycle;
\fill [fill=lightgray] (12,0)--(12,6)--(16,6)--(16,0)--cycle;

\draw[very thick] (0,0)--(0,6);

\draw[thick, densely dotted] (4,0)--(4,6);
\draw[very thick] (6,0)--(6,6);

\draw[thick, densely dotted] (10,0)--(10,6);
\draw[very thick] (12,0)--(12,6);

\draw[thick, densely dotted] (16,0)--(16,6);
\draw[very thick] (18,0)--(18,6);

\draw[densely dotted] (20,0)--(20,6);
\draw[solid] (18,0)--(20,6);
\draw[solid] (22,0)--(20,6);

\draw[densely dotted] (22,0)--(22,6);
\draw[solid] (22,0)--(24,6);
\draw[very thick] (24,0)--(24,6);

\draw[solid] (0,6)--(24,6);
\draw[solid] (0,0)--(24,0);
\end{tikzpicture}
\caption{An illustration of sequences:
$(\mb{u}_i)_{i \in [n]}$ and $(\mb{\wt u}_i)_{i \in [n]}$ are aligned on the upper row and on the lower row respectively (above), and similarly $(\mb{v}_i)_{i \in [N-n]}$ and $(\mb{\wt v}_i)_{i \in [N-n]}$ are aligned (below).}
\label{Fig:coupling_matrix}
\end{figure}

Let us define two sequences $(\mb{u}_i)_{i \in [n]}$ and $(\mb{\wt u}_i)_{i \in [n]}$ 
of elements in $S \times S$,
aligning $(u, v)\in S\times S$ depending on $\xi^{\sigma_0}_{u, v}(\sigma)$ and $\xi^{\sigma_0}_{u, v}(\wt \sigma)$.
See also Figure \ref{Fig:coupling_matrix}.
First we divide $[n]$ into intervals $O_u$ of size $\xi_u(\sigma_0)$ for $u \in S$.
Second each $O_u$ is decomposed into intervals $O_{u,v}$ for $v \in S$ of size $m_{u,v}$
and $O_u^+$ of size $\xi_u(\sigma_0)-\sum_{v \in S}m_{u,v}$.
For $u, v \in S$,
let us define $\mb{u}_i=(u, v)$ for $i \in O_{u,v}$,
and $(u, w)$ for $i \in O_u^+$ and for some $w \in S$
so that for each $u \in S$,
the total number of $(u, v)$ in $O_u$ is $\xi_{u,v}^{\sigma_0}(\sigma)$.
Analogously, we define $\mb{\wt u}_i$ for $i \in O_u$ by using $\xi^{\sigma_0}_{u, v}(\wt \sigma)$.

By \eqref{Eq:Lem:coupling_matrix:1}, we have $\#O_{u, v}=m_{u, v} \ge \ceil{n/(2Q^2)}$.
Further each $O_{u,v}$ is divided into $O_{u,v}^o$ of size $\ceil{n/(2Q^2)}$ and $O_{u,v}^r$ of size $m_{u,v}-\ceil{n/(2Q^2)}$.
Letting $\wt O_u^+=O_u^+$,
we decompose $O_u^+$ into $O_{u, v}^+$ of size equal to the number of $i \in O_u^+$
with $\mb{u}_i=(u, v)$.
Analogously,
we decompose $\wt O_u^+$ into $\wt O_{u,v}^+$ of size equal to the number of $i \in \wt O_u^+$ with $\mb{\wt u}_i=(u, v)$.
Note that either $O_{u,v}^+$ or $\wt O_{u,v}^+$ (possibly both) is empty.

Let us define two sequences $(\mb{v}_i)_{i \in [N-n]}$ and $(\mb{\wt v}_i)_{i \in [N-n]}$ of elements in $S$,
where the former contains $\kappa-\xi_v(\sigma)$ copies of $v$ and the latter contains $\kappa-\xi_v(\wt \sigma)$ copies of $v$ for $v \in S$.
These sequences indicate the number of ``vacancies'' for each $v \in S$.
We align the sequences so that for a set $V_v$ of indices of size $\kappa-M_v$,
we have $\mb{v}_i=\mb{\wt v}_i=v$ for $i \in V_v$.
In $[N-n]\setminus \bigsqcup_{v \in S}V_v$,
let $V_v^+$ be a subset of size $M_v-\xi_v(\sigma)$ such that $\mb{v}_i=v$ for $i \in V_v^+$,
and $\wt V_v^+$ be a subset of size $M_v-\xi_v(\wt \sigma)$ such that $\mb{\wt v}_i=v$ for $i \in \wt V_v^+$. 
Note that at least $V_v^+$ or $\wt V_v^+$ is empty for each $v \in S$,
and $\bigsqcup_{v \in S}V_v^+=\bigsqcup_{v \in S}\wt V_v^+$.
Further if $i \in V_v^+$,
then $\mb{v}_i=v$ and $\mb{\wt v}_i=w$ for some $w \neq v$.

By \eqref{Eq:Lem:coupling_matrix:1},
for each $v \in S$,
we have $\#V_v=\kappa-M_v \ge \ceil{(N-n)/(2Q)}$.
Further each $V_v$ is divided into $V_v^o$ of size $\ceil{(N-n)/(2Q)}$ and $V_v^r:=V_v\setminus V_v^o$.
Let
\[
O^+:=\bigsqcup_{u \in S}O_u^+ \quad \text{and} \quad V^+:=\bigsqcup_{v \in S}V_v^+.
\]
Note that $\mb{D}^{\sigma_0}(\sigma, \wt \sigma)=\#O^+$ and $D^{\sigma_0}_u(\sigma, \wt \sigma)=\#O_u^+$ for each $u\in S$,
and further that $D(\sigma, \wt \sigma)=\#V^+$.
By the triangle inequality, 
it holds that $D(\sigma, \wt \sigma)\le \mb{D}^{\sigma_0}(\sigma, \wt \sigma)$,
and thus 
\begin{equation}\label{Eq:V+O+}
\#V^+\le \#O^+.
\end{equation}

The transition is defined as follows.
The state is kept with probability $1/(1+N-n)$,
or the following is performed with probability $(N-n)/(1+N-n)$.
We say that {\bf the move $(\mb{u}, \mb{v})=(u_\ast, v_\ast, v_\ast')\in (S\times S)\times S$
is applied to $\mb{\xi}(\sigma)=\{\xi_{u, v}(\sigma)\}_{u,v \in S}$}
when $\xi_{u_\ast, v_\ast}(\sigma)$ is decremented by $1$ and then $\xi_{u_\ast, v_\ast'}(\sigma)$
is incremented by $1$.
(If $v_\ast=v_\ast'$,
then $\xi_{u_\ast, v_\ast}(\sigma)$ does not change.)
Choosing $i \in [n]$ and $j \in [N-n]$ uniformly at random respectively and independently,
we will apply the move $(\mb{u}_i, \mb{v}_j)$ to $\mb{\xi}(\sigma)$ and some move to $\mb{\xi}(\wt \sigma)$ respectively.

\subsubsection*{\upshape{(i) {\bf Case} $(i, j)\in O_{u_\ast, v_\ast}\times V_{v_\ast'}$}}
First let us consider the case when $(i, j)\notin O^o_{u_\ast, v_\ast}\times V^o_{v_\ast'}$.
The moves $(\mb{u}_i, \mb{v}_j)$ and $(\mb{\wt u}_i, \mb{\wt v}_j)$ are applied to $\mb{\xi}(\sigma)$ and $\mb{\xi}(\wt \sigma)$ respectively.
The move does not change $\mb{D}^{\sigma_0}(\sigma, \wt \sigma)$.
Second let us consider the case when $(i, j)\in O^o_{u_\ast, v_\ast}\times V^o_{v_\ast'}$.
If $u_\ast=u$ and $v_\ast, v_\ast' \in \{v_0, v_1\}$ hold,
then the moves are applied as in the following,
otherwise the move $(\mb{u}_i, \mb{v}_j)$ and $(\mb{\wt u}_i, \mb{\wt v}_j)$ are applied to $\mb{\xi}(\sigma)$ and $\mb{\xi}(\wt \sigma)$ respectively:

If $(v_\ast, v_\ast')=(v_k, v_k)$ for $k=0, 1$,
then the move $(u, v_k, v_k)$ is applied to $\mb{\xi}(\sigma)$ and the move $(u, v_k, v_{1-k})$ is applied to $\mb{\xi}(\wt \sigma)$.
Further if $(v_\ast, v_\ast')=(v_k, v_{1-k})$ for $k=0, 1$,
then the move $(u, v_k, v_{1-k})$ is applied to $\mb{\xi}(\sigma)$ and the move $(u, v_k, v_k)$ is applied to $\mb{\xi}(\wt \sigma)$.
This move decreases or increases $\mb{D}^{\sigma_0}(\sigma, \wt \sigma)$ by $1$,
and it happens with equal probability
\[
\frac{2}{n(N-n)}\ceil{\frac{n}{2Q^2}}\ceil{\frac{N-n}{2Q}}\frac{N-n}{1+N-n}\ge \frac{1}{2Q^3}\(1-\frac{1}{1+N-n}\)\ge \frac{1}{4Q^3}.
\]
\subsubsection*{\upshape{(ii) {\bf Case} $(i, j)\in O^+_{u_\ast, v_\ast}\times V_{v_\ast'}\sqcup O_{u_\ast, v_\ast}\times V_{v_\ast'}^+$}}
First suppose that $(i, j) \in O^+_{u_\ast, v_\ast}\times V_{v_\ast'}$,
in which case $\mb{u}_i=(u_\ast, v_\ast)$ and $\mb{v}_j=v_\ast'$,
and $\mb{\wt u}_i=(u_\ast, \wt v_\ast)$ and $\mb{\wt v}_j=v_\ast'$,
where $\wt v_\ast \neq v_\ast$.
The moves $(u_\ast, v_\ast, v_\ast')$ and $(u_\ast, \wt v_\ast, v_\ast')$ are applied to $\mb{\xi}(\sigma)$ and $\mb{\xi}(\wt \sigma)$ respectively.
This decreases $\mb{D}^{\sigma_0}(\sigma, \wt \sigma)$ by $1$,
and it happens with probability
\[
\frac{\#O^+}{n}\frac{N-n-\#V^+}{N-n}\frac{N-n}{1+N-n}.
\]

Next suppose that $(i, j) \in O_{u_\ast, v_\ast}\times V_{v_\ast'}^+$,
in which case $\mb{u}_i=(u_\ast, v_\ast)$ and $\mb{v}_j=v_\ast'$, and $\mb{\wt u}_i=(u_\ast, v_\ast)$ and $\mb{\wt v}_j=\wt v_\ast'$ for some $\wt v_\ast'\neq v_\ast'$.
The moves $(u_\ast, v_\ast, v_\ast')$ and $(u_\ast, v_\ast, \wt v_\ast')$ are applied to $\mb{\xi}(\sigma)$ and $\mb{\xi}(\wt \sigma)$ respectively.
In this case, the move possibly increases $\mb{D}^{\sigma_0}(\sigma, \wt \sigma)$ by $1$.
This happens with probability
\[
\frac{n-\#O^+}{n}\frac{\#V^+}{N-n}\frac{N-n}{1+N-n}.
\]
The overall effect to $\mb{D}^{\sigma_0}(\sigma, \wt \sigma)$ on average up to factor $(N-n)/(1+N-n)$ is at most
\begin{align*}
-\frac{\#O^+}{n}\frac{N-n-\#V^+}{N-n}+\frac{n-\#O^+}{n}\frac{\#V^+}{N-n}
=\frac{-\#O^+(N-n)+n\#V^+}{n(N-n)}\le 
\frac{-\#O^+(N-2n)}{n(N-n)},
\end{align*}
where we have used $\#V^+\le \#O^+$ by \eqref{Eq:V+O+} in the last inequality.
The last term is at most $0$ by the assumption that $N \ge 2n$, following since $n=\floor{\rho N}\le N/2$ for $\rho \in (0, 1/2]$.
\subsubsection*{\upshape{(iii) {\bf Case} $(i, j) \in O_{u_\ast, v_\ast}^+\times V_{v_\ast'}^+$}}
In this case $\mb{u}_i=(u_\ast, v_\ast)$ and $\mb{\wt u}_i=(u_\ast, \wt v_\ast)$ for some $\wt v_\ast \neq v_\ast$,
and $\mb{v}_j=v_\ast'$ and $\mb{\wt v}_j=\wt v_\ast'$ for some $\wt v_\ast'\neq v_\ast'$.
The moves $(u_\ast, v_\ast, v_\ast')$ and $(u_\ast, \wt v_\ast, \wt v_\ast')$ are applied to $\mb{\xi}(\sigma)$ and $\mb{\xi}(\wt \sigma)$ respectively.
In this case $\mb{D}^{\sigma_0}(\sigma, \wt \sigma)$ does not increase,
more precisely, the possible increments are $-2$, $-1$ or $0$.

Summarizing all the cases concludes the claim.
\qed

\begin{theorem}\label{Thm:upper}
Let $Q$ be a fixed integer at least two, $\rho$ be a fixed real in $(0,1/2]$.
For a positive integer $\kappa$,
let $N=\kappa Q$ and $n=\floor{\rho N}$.
For the chain $\{\sigma(t)\}_{t\in \Z_+}$ on $\Sc_n$,
letting $t_n:=\ceil{(1/2)n\log n}$,
we have
\[
\lim_{\alpha \to \infty}\limsup_{n\to \infty}D_\TV(t_n +\ceil{\alpha n})=0.
\]
\end{theorem}

\proof
For each $\sigma(0) \in \Sc_n$,
let us define a coupling between the chains $\{\sigma(t)\}_{t\in \Z_+}$ starting from $\sigma(0)$ and $\{\wt \sigma(t)\}_{t\in \Z_+}$ with the initial distribution $\mu_n$ respectively.
For $\alpha>0$,
let $T_{\alpha, N}:=\ceil{\alpha N}$.
We run two chains independently up to time $T_{\alpha, N}$ and after that we consider
two independent load matrix chains up to time $t_n+T_{\alpha, N}$, and then apply to the coupling in Lemma \ref{Lem:coupling_matrix}.

Let us fix a constant $\delta>0$ satisfying that
$0<\max\big\{\delta, c_\star\sqrt{\delta}\big\}\le \rho/(8Q)$,
where $c_\star$ is the constant depending only on $\rho$ and $Q$ in Lemma \ref{Lem:load_matrix_good}.
By the first claim of Lemma \ref{Lem:burning-in},
for all large enough $\alpha$ and all large enough $N$,
\[
\Pb_{\sigma(0)}\big(\sigma(T_{\alpha, N}) \notin \Sc(\delta)\big) \le e^{-c_\delta N}.
\]
Let $\sigma_0:=\sigma(T_{\alpha, N})$.
By Lemma \ref{Lem:TV_stationary},
for $t \in \Z_+$,
\begin{equation}\label{Eq:TV}
\|\Pb_{\sigma_0}\big(\sigma(t)\in \cdot\,\big)-\mu_n\|_\TV=\|\Pb_{\sigma_0}\big(\xi^{\sigma_0}(t)\in \cdot\,\big)-\wbar \mu_n^{\sigma_0}\|_\TV.
\end{equation}
Let us consider the case when $\sigma_0 \in \Sc(\delta)$, 
if otherwise we run the chains independently.
Letting
\[
\tau_\coup^{\sigma_0}:=\inf\big\{t \ge 0 \ : \ \xi^{\sigma_0}(\sigma(t))=\xi^{\sigma_0}(\wt \sigma(t))\big\},
\]
we will show that under this coupling the probability that $\tau_\coup^{\sigma_0}>t_n+\alpha n$ tends to $0$ as $n \to \infty$ and then as $\alpha \to \infty$.

By Lemma \ref{Lem:load_matrix_square_root},
we have $C$ such that for all large enough $n$,
for all $\sigma_\ast \in \Sc_n$ and $R>1$,
\[
\Pb_{\sigma_\ast}\(\sigma(t_n) \notin \Sc\(\sigma_0, \frac{R}{\sqrt{n}}\)\) \le \frac{CQ^2}{R}.
\]
This implies that for all large enough $n$ and for $R>1$,
\begin{equation}\label{Eq:Thm:upper:square_root}
\Pb_{\sigma_0}\(\sigma(t_n) \notin \Sc\(\sigma_0, \frac{R}{\sqrt{n}}\)\) \le \frac{CQ^2}{R}, \quad
\Pb_{\mu_n}\(\wt \sigma(t_n) \notin \Sc\(\sigma_0, \frac{R}{\sqrt{n}}\)\) \le \frac{CQ^2}{R}.
\end{equation}

For each $R>1$ for all large enough $n$,
we have $R/\sqrt{n}\le \delta$.
For $t \in \Z_+$, 
let 
\[
\mb{G}_t:=\Big\{\sigma(t) \in \Sc(2\delta)\cap \Sc\big(\sigma_0, c_\star \sqrt{\delta}\big)\Big\}
\quad \text{and} \quad 
\mb{\wt G}_t:=\Big\{\wt \sigma(t) \in \Sc(2\delta)\cap \Sc\big(\sigma_0, c_\star \sqrt{\delta}\big)\Big\}.
\]
By the second claim of Lemma \ref{Lem:burning-in}, and by Lemma \ref{Lem:load_matrix_good},
for $\sigma_0  \in \Sc(\delta)$ and $\sigma_\ast, \wt \sigma_\ast \in \Sc(\sigma_0, \delta)$,
then
\begin{equation}\label{Eq:GwtG}
\Pb_{\sigma_\ast}\Big(\bigcup_{0\le t \le n^2}\mb{G}_t^{\sf c}\Big)\le e^{-c_\delta n}
\quad \text{and} \quad
\Pb_{\wt \sigma_\ast}\Big(\bigcup_{0\le t\le n^2}\mb{\wt G}_t^{\sf c}\Big)\le e^{-c_\delta n}.
\end{equation}
Let us assume that $\mb{G}_t\cap \mb{\wt G}_t$ holds,
if otherwise we run the chains independently.
If $\mb{D}^{\sigma_0}_t=0$,
then we run them with common transitions after $t$.
By Lemma \ref{Lem:coupling_matrix},
as long as $\mb{D}^{\sigma_0}_t>0$ and $\mb{G}_t\cap \mb{\wt G}_t$ holds,
then the variance $\Var\(\mb{D}^{\sigma_0}_{t+1}-\mb{D}^{\sigma_0}_t \mid \sigma(t), \wt \sigma(t)\)$ is at least $1/(2Q^3)$.
Thus, for $\alpha>0$ and $R>1$ for all large enough $n$ the following holds by \eqref{Eq:GwtG},
\begin{align*}
\Pb_{\sigma_\ast,\wt \sigma_\ast}\big(\tau_\coup^{\sigma_0}>\alpha n\big)
&\le \Pb_{\sigma_\ast, \wt \sigma_\ast}\Big(\big\{\tau_\coup^{\sigma_0}>\alpha n\big\}\cap \bigcap_{0\le t\le \ceil{\alpha n}}\mb{G}_t\cap \mb{\wt G}_t\Big)+2e^{-c_\delta n}\\
&\le \frac{4\mb{D}^{\sigma_0}_0}{\sqrt{(2Q^3)^{-1}\alpha n}}+2e^{-c_\delta n}
\le \frac{8Q^{3/2}R}{\sqrt{\alpha}}+2e^{-c_\delta n}\le \frac{10Q^{3/2}R}{\sqrt{\alpha}}.
\end{align*}
In the above we have used \cite[Proposition 17.20]{LP} in the second inequality,
and $\mb{D}^{\sigma_0}_0\le \sum_{u \in S}(R/\sqrt{n})\xi_u(\sigma_0)= R\sqrt{n}$ in the third inequality.
Thus,
if $\sigma_0\in \Sc(\delta)$ and $\sigma(t_n), \wt \sigma(t_n) \in \Sc(\sigma_0, R/\sqrt{n})$,
then for $\alpha>0$,
\[
\Pb_{\sigma(t_n), \wt \sigma(t_n)}\(\tau_\coup^{\sigma_0}>\alpha n\)\le \frac{10Q^{3/2}R}{\sqrt{\alpha}}.
\]
Furthermore, since $\sigma_0=\sigma(T_{\alpha, N}) \in \Sc(\delta)$,
by the Markov property and by \eqref{Eq:Thm:upper:square_root},
\begin{align*}
&\left\|\Pb_{\sigma_0}\big(\xi^{\sigma_0}(t_n+\ceil{\alpha n})\in \cdot\,\big)-\wbar \mu_n^{\sigma_0}\right\|_\TV\le 
\Pb_{\sigma_0, \mu_n}\big(\tau_\coup^{\sigma_0}>t_n+\alpha n\big)\\
&\le \Eb\Big[\Pb_{\sigma(t_n), \wt \sigma(t_n)}\big(\tau_\coup^{\sigma_0}>\alpha n\big)\1_{\{\sigma(t_n), \wt \sigma(t_n)\in \Sc(\sigma_0, R/\sqrt{n})\}}\Big]+\frac{2CQ^2}{R}
\le \frac{10Q^{3/2}R}{\sqrt{\alpha}}+\frac{2CQ^2}{R}.
\end{align*}
Therefore it holds that
\begin{align*}
&\left\|\Pb_{\sigma(0)}\big(\sigma(T_{\alpha, N}+t_n+\ceil{\alpha n})\in \cdot\,\big)-\mu_n\right\|_\TV\\
&\le \Eb\Big[\big\|\Pb_{\sigma_0}\big(\sigma(t_n+\ceil{\alpha n})\in \cdot\,\big)-\mu_n\big\|_\TV\1_{\{\sigma_0\in \Sc(\delta)\}}\Big]+e^{-c_\delta n}
\le \frac{10Q^{3/2}R}{\sqrt{\alpha}}+\frac{2CQ^2}{R}+e^{-c_\delta n}.
\end{align*}
In the above we have used the Markov property and \eqref{Eq:TV}.
Taking maximum over the initial state $\sigma(0)$ yields for $\alpha>0$ and $R>1$ and for all large enough $n$,
\[
D_{\TV}(t_n+T_{\alpha, N}+\ceil{\alpha n})\le
 \frac{10Q^{3/2}R}{\sqrt{\alpha}}+\frac{2CQ^2}{R}+e^{-c_\delta n}. 
\]
Letting $n \to \infty$ and then $\alpha\to \infty$ and further $R \to \infty$ concludes the claim.
\qed

\begin{theorem}\label{Thm:lower}
Let $Q$ be a fixed integer at least two,
$\rho$ be a fixed real in $(0, 1)$.
For a positive integer $\kappa$,
let $N=\kappa Q$ and $n=\floor{\rho N}$.
For the chain $\{\sigma(t)\}_{t\in \Z_+}$ on $\Sc_n$,
letting $t_n:=\ceil{(1/2)n\log n}$,
we have
\[
\lim_{\alpha\to \infty}\liminf_{n \to \infty}D_\TV\(t_n-\floor{\alpha n}\)=1.
\]
\end{theorem}

\proof
Fix an initial state $\sigma_0\in \Sc_n$ satisfying that for all $u \in S$,
\[
\left|\frac{\xi_u(\sigma_0)}{n}-\frac{1}{Q}\right|\le \frac{1}{n}.
\]
Note that such $\sigma_0$ exists for each $n>0$.
By \eqref{Eq:load_matrix} and the expectations,
it holds that
\begin{align*}
\xi^{\sigma_0}_{u,v}(t+1)-\Eb \xi^{\sigma_0}_{u,v}(t+1)
&=\(\frac{\kappa-\xi_v(t)}{(1-\rho)N}-\Eb\frac{\kappa-\xi_v(t)}{(1-\rho)N}\)\frac{\xi_u(\sigma_0)}{n}\\
&+\(1-\frac{1}{n}\)\(\xi^{\sigma_0}_{u,v}(t)-\Eb \xi^{\sigma_0}_{u,v}(t)\)+\wh M_{u,v,t+1}-\Eb\wh M_{u,v,t+1}.
\end{align*}
In the above $|\wh M_{u,v,t+1}|\le C$ and $|\Eb[\wh M_{u,v,t+1}\mid \Fc_t]|\le C/n$ almost surely.
Squaring both sides yields the following:
\begin{align*}
&\(\xi^{\sigma_0}_{u,v}(t+1)-\Eb \xi^{\sigma_0}_{u,v}(t+1)\)^2\\
&=\(-\frac{\xi_v(t)}{(1-\rho)N}+\Eb\frac{\xi_v(t)}{(1-\rho)N}\)^2\(\frac{\xi_u(\sigma_0)}{n}\)^2+\(1-\frac{1}{n}\)^2\(\xi^{\sigma_0}_{u,v}(t)-\Eb \xi^{\sigma_0}_{u,v}(t)\)^2\\
&+2\(-\frac{\xi_v(t)}{(1-\rho)N}+\Eb\frac{\xi_v(t)}{(1-\rho)N}\)\frac{\xi_u(\sigma_0)}{n}\(1-\frac{1}{n}\)\(\xi^{\sigma_0}_{u,v}(t)-\Eb \xi^{\sigma_0}_{u,v}(t)\)+\wt M_{u,v,t+1}.
\end{align*}
In the above $|\wt M_{u,v,t+1}|\le C'n$ and $|\Eb[\wt M_{u,v,t+1}\mid \Fc_t]|\le C'$ almost surely.
Summing over $u \in S$ yields the following:
\begin{align*}
&\sum_{u\in S}\(\xi^{\sigma_0}_{u,v}(t+1)-\Eb \xi^{\sigma_0}_{u,v}(t+1)\)^2\\
&=\(\frac{\xi_v(t)}{(1-\rho)N}-\Eb\frac{\xi_v(t)}{(1-\rho)N}\)^2\(\frac{1}{Q}+O\(\frac{1}{n}\)\)
+\(1-\frac{1}{n}\)^2\sum_{u\in S}\(\xi^{\sigma_0}_{u,v}(t)-\Eb \xi^{\sigma_0}_{u,v}(t)\)^2\\
&-2\(\frac{\xi_v(t)}{(1-\rho)N}-\Eb\frac{\xi_v(t)}{(1-\rho)N}\)\(1-\frac{1}{n}\)\frac{1}{Q}\sum_{u \in S}\(\xi^{\sigma_0}_{u,v}(t)-\Eb \xi^{\sigma_0}_{u,v}(t)\)+O(1)+\sum_{u \in S}\wt M_{u,v,t+1}.
\end{align*}
Noting that
$\sum_{u\in S}\(\xi^{\sigma_0}_{u,v}(t)-\Eb \xi^{\sigma_0}_{u,v}(t)\)=\xi_v(t)-\Eb \xi_v(t)$,
we have by the expectations,
\begin{align*}
\sum_{u \in S}\Var\(\xi^{\sigma_0}_{u,v}(t+1)\)
&=\(1-\frac{1}{n}\)^2\sum_{u \in S}\Var\(\xi^{\sigma_0}_{u,v}(t)\)+\frac{1}{Q(1-\rho)^2 N^2}\Var\(\xi_v(t)\)\\
&-2\(1-\frac{1}{n}\)\frac{1}{Q(1-\rho)N}\Var\(\xi_v(t)\)+O(1).
\end{align*}
We note that for all large enough $n$,
\[
\frac{1}{Q(1-\rho)^2N^2}<2\(1-\frac{1}{n}\)\frac{1}{Q(1-\rho)N}.
\]
Therefore it holds that
\[
\sum_{u \in S}\Var\(\xi^{\sigma_0}_{u,v}(t+1)\)\le \(1-\frac{1}{n}\)^2\sum_{u \in S}\Var\(\xi^{\sigma_0}_{u,v}(t)\)+C.
\]
Iterating in $t \in \Z_+$ shows that since $\Var\(\xi^{\sigma_0}_{u,v}(0)\)=0$,
\[
\sum_{u\in S}\Var\(\xi^{\sigma_0}_{u,v}(t)\)
\le 
\(1-\frac{1}{n}\)^{2t}\sum_{u\in S}\Var\(\xi^{\sigma_0}_{u,v}(0)\)+Cn=Cn.
\]
In the normalized form $\wbar \xi^{\sigma_0}_{u,v}(t)$,
for each $v \in S$, for all $t \in \Z_+$ and for all large enough $n$,
\[
\sum_{u \in S}\Var\(\wbar \xi^{\sigma_0}_{u,v}(t)\)\le \frac{C'}{n}.
\]
In particular, letting $t \to \infty$ yields
\[
\sum_{u \in S}\Var_{\mu_n}\(\wbar \xi^{\sigma_0}_{u,v}\)\le \frac{C'}{n}.
\]
Let us fix an arbitrary $u \in S$, and focus on the case when $v=u$.
By the assumption on the initial state $\sigma_0$ and by \eqref{Eq:load_eq},
for all large enough $n$ and all $t\in \Z_+$,
\[
\left|\Eb \frac{\xi_u(t)}{n}-\frac{1}{Q}\right| \le \frac{C}{n}.
\]
This combined with \eqref{Eq:load_matrix_eq} implies that for all $t \in \Z_+$,
\[
\Eb\wbar \xi^{\sigma_0}_{u,u}(t+1)-\frac{1}{Q}=
\(1-\frac{1}{n}\)\(\Eb \wbar \xi^{\sigma_0}_{u,u}(t)-\frac{1}{Q}\)+O\(\frac{1}{n^2}\).
\]
Thus by iteration in $t \in \Z_+$ yields
\[
\Eb\wbar \xi^{\sigma_0}_{u,u}(t)=\frac{1}{Q}+\(1-\frac{1}{n}\)^t\(1-\frac{1}{Q}\)+O\(\frac{1}{n}\).
\]
In the above we have used $\wbar \xi^{\sigma_0}_{u,u}(0)=1$ by the assumption on $\sigma_0$.
For $t_{\alpha, n}:=t_n-\floor{\alpha n}$ for $\alpha>0$,
we have
\[
\Eb\wbar \xi^{\sigma_0}_{u,u}(t_{\alpha,n})=\frac{1}{Q}+(1+o(1))\frac{e^\alpha}{\sqrt{n}}\(1-\frac{1}{Q}\).
\]
The estimate on the variance implies that $\Var\(\wbar \xi^{\sigma_0}_{u,u}(t_{\alpha, n})\) \le C'/n$ and thus for all $R>0$,
\[
\Pb\(\wbar \xi^{\sigma_0}_{u,u}(t_{\alpha, n}) \ge \frac{1}{Q}+(1+o(1))\(1-\frac{1}{Q}\)\frac{e^\alpha}{\sqrt{n}}-\frac{R}{\sqrt{n}}\)=1-O\(\frac{1}{R^2}\).
\]
Further under the stationary distribution $\mu_n$, 
since $\Var_{\mu_n}\(\wbar \xi^{\sigma_0}_{u,u}\)\le C'/n$,
we have
\[
\mu_n\(\wbar \xi^{\sigma_0}_{u,u} \ge \frac{1}{Q}+\frac{R}{\sqrt{n}}\) =O\(\frac{1}{R^2}\).
\]
Summarizing the above estimates concludes that for all $\alpha>0$ and for all large enough $n$,
\[
D_\TV(t_{\alpha, n})=1-O\(e^{-2\alpha}\).
\]
Letting $n \to \infty$ and then $\alpha\to \infty$ shows the claim.
\qed

\proof[Proof of Theorem \ref{Thm:Gibbs}]
The claim follows from Theorems \ref{Thm:upper} and \ref{Thm:lower}.
\qed

\appendix

\section{Uniform sampling on M-convex sets}\label{Sec:appendix}

Let $Q$ and $n$ be positive integers.
An {\bf M-convex set} $B$ is a subset of
\[
\left\{(x_1, \dots, x_Q) \in \Z^Q \ : \ \sum_{i=1}^Q x_i=n\right\}
\]
such that the following holds:
For all $x, y \in B$ with $x \neq y$,
for every $i\in [Q]$ with $x_i>y_i$ 
there exists $j \in [Q]$ with $y_j>x_j$,
satisfying that $x-\mb{e}_i+\mb{e}_j \in B$ and $y-\mb{e}_j+\mb{e}_i \in B$.
For example, 
for all positive integers $Q$ and $n$,
the following set is M-convex:
\[
B^{(Q)}(n):=\left\{(x_1, \dots, x_Q) \in \Z_+^Q \ : \ \sum_{i=1}^Q x_i=n\right\}.
\]
This set coincides with $\Xi_{n, S, \kappa}$ if $n\le \kappa$;
the state space of load profile chains
for positive integers $\kappa$, $Q$ and $n$
(for the notation, see Section \ref{Sec:introduction}).
Note that for $n \le Q$,
\[
\Ic:=\left\{(x_1, \dots, x_Q) \in \{0, 1\}^Q \ : \ \sum_{i=1}^Q x_i \le n\right\}
\]
defines a (uniform) matroid $\Mcc=\([Q], \Ic\)$ \cite[Example 1.2.7, Chapter 2]{Oxley}.
Moreover,
\[
B_{\Ic}:=\left\{(x_1, \dots, x_Q) \in \{0, 1\}^Q \ : \ \sum_{i=1}^Qx_i =n\right\}
\]
forms the set of bases in $\Mcc$.
The notion of M-convex set generalizes the set of bases in a matroid.
See \cite{Murota} for background.

Every M-convex set $B$
admits a natural Markov chain with the uniform stationary distribution.
The transition probability is defined as follows:
Let $x=(x_1, \dots, x_Q) \in B$.
An $i\in [Q]$ is chosen uniformly at random,
and then independently a $j$ is chosen uniformly at random from
\[
\Big\{j\in [Q] \ : \ x-\mb{e}_i+\mb{e}_j \in B\Big\}.
\]
Note that the above set is not empty since it contains $i$.
Letting $y:=x-\mb{e}_i+\mb{e}_j$ in $B$
defines the transition from $x$ to $y$.
The chain is irreducible since $B$ is M-convex,
and is aperiodic since it has a holding probability at each state.
Further it is symmetric.
Indeed,
for $x, y \in B$ with $y=x-\mb{e}_i+\mb{e}_j$,
the transition from $x$ to $y$ occurs with probability
\[
\frac{1}{Q\#\{j \in [Q] \ : \ x-\mb{e}_i+\mb{e}_j \in B\}}.
\]
This coincides with the transition probability from $y$ to $x$.
In the case when $x_i=0$,
then $j=i$ almost surely and $y=x$.
Therefore the chain is reversible with respect to the uniform distribution,
which is the unique stationary distribution.
Let us restrict ourselves to the special case when
\[
B=B^{(Q)}(n).
\]
A particular feature of this case is that given a state $x$,
if the first index $i$ satisfies that $x_i>0$,
then the second index $j$ is always chosen uniformly from the whole $[Q]$.
If we identify $B^{(Q)}(n)$ with $\Xi_{n, S, \kappa}$ for $n \le \kappa$,
then
the uniform distribution ${\rm Unif}_n$ on $\Xi_{n, S, \kappa}$ is obtained as follows:
For an arbitrary real $\lambda>0$,
\[
{\rm Unif}_n(\xi)=\frac{1}{Z_\lambda}\prod_{v \in S}\frac{\lambda^{\xi_v}}{1+\lambda+\cdots+\lambda^\kappa} \quad \text{where $Z_\lambda:= \frac{\lambda^n}{(1+\lambda+\cdots+\lambda^\kappa)^Q}$}\cdot\# \Xi_{n, S, \kappa}.
\]

We consider the lazy version of the chain (i.e.\ the original process is delayed with $1/2$-holding probability in each step) for an upper bound on the mixing time in total variation.

\begin{proposition}\label{Prop:uniform}
For each integer $Q$ at least two and for each real $\e\in (0, 1)$,
there exists a constant $C_\e>0$ such that for all $n>0$,
\[
T^\mix(\e)\le C_\e n^2.
\]
\end{proposition}

\proof
Let $\{X_t\}_{t \in \Z_+}$ and $\{Y_t\}_{t\in \Z_+}$ be the lazy Markov chains with the common transition probability,
and we construct a coupling between them in the following.
Suppose that $X_t=x$ and $Y_t=y$.
If $x=y$,
then we run the chains with common transitions.
Otherwise if $x\neq y$,
then we flip an independent fair coin $\chi \in \{{\sf head}, {\sf tail}\}$ and independently choose $i \in [Q]$ uniformly at random.
Let 
\[
\left\{{\text{$x=0$ or $y=0$}}\right\}:=\left\{i \in [Q] \ : \ x_i\cdot y_i=0\right\},
\quad
\left\{x>y>0\right\}:=\left\{i \in [Q] \ : \ x_i>y_i>0\right\}.
\]
Similarly, we define
\[
\left\{x=y>0\right\}:=\left\{i \in [Q] \ : \ x_i=y_i>0\right\},
\quad
\left\{y>x>0\right\}:=\left\{i \in [Q] \ : \ y_i>x_i>0\right\}.
\]
Analogously we write $\left\{x>y=0\right\}$ and $\left\{y>x=0\right\}$.
It holds that
\[
[Q]=\left\{x>y>0\right\}\bigsqcup\left\{y>x>0\right\}\bigsqcup\left\{\text{$x=0$ or $y=0$}\right\}\bigsqcup\left\{x=y>0\right\}.
\]
In the above some of the sets are possibly empty.
\subsubsection*{\upshape{\bf (1) Case} $i\in \left\{x>y>0\right\}\bigsqcup\left\{y>x>0\right\}$}
The second index $j$ is chosen from $[Q]$ uniformly at random both for $X_t$ and for $Y_t$.
\subsubsection*{\upshape{\bf (1a) Case} $j \in \left\{x>y>0\right\}\bigsqcup\left\{y>x>0\right\}$}
If $\chi={\sf head}$,
then we run $X_t$ to the next with the chosen $i$ and $j$, and keep $Y_t$.
Namely,
\[
X_{t+1}:=X_t-\mb{e}_i+\mb{e}_j \quad \text{and} \quad Y_{t+1}:=y.
\]
Further,
if $\chi={\sf tail}$,
then we change the roles of $X_t$ and $Y_t$ and run $Y_t$ to the chosen $i$ and $j$,
and keep $X_t$.
\subsubsection*{\upshape{\bf (1b) Case} $j \notin \left\{x>y>0\right\}\bigsqcup\left\{y>x>0\right\}$} If $\chi={\sf head}$,
then we move $X_t$ and $Y_t$ together with the chosen $i$ and $j$.
Namely,
\[
X_{t+1}:=X_t-\mb{e}_i+\mb{e}_j \quad \text{and} \quad Y_{t+1}:=Y_t-\mb{e}_i+\mb{e}_j,
\]
and if $\chi={\sf tail}$,
then we keep both, i.e., $X_{t+1}:=x$ and $Y_{t+1}:=y$.
\subsubsection*{\upshape{\bf (2) Case} $i \in \left\{\text{$x=0$ or $y=0$}\right\}$}
If $\chi={\sf head}$,
then we run $X_t$ to the next with the chosen $i$ and a chosen $j$ according to the rule,
and keep $Y_t$.
Namely,
\[
X_{t+1}:=X_t-\mb{e}_i+\mb{e}_j \quad \text{and} \quad Y_{t+1}:=y.
\]
Further, if $\chi={\sf tail}$,
then we change the roles of $X_t$ and $Y_t$ and run $Y_t$ to the next with the chosen $i$ and a chosen $j$ according to the rule, and keep $X_t$.
\subsubsection*{\upshape{\bf (3) Case} $i \in \left\{x=y>0\right\}$}
The second index $j$ is chosen from $[Q]$ uniformly at random both for $X_t$ and for $Y_t$.
If $\chi={\sf head}$,
then we move $X_t$ and $Y_t$ together with the chosen $i$ and $j$.
Namely,
\[
X_{t+1}:=X_t-\mb{e}_i+\mb{e}_j
\quad \text{and} \quad
Y_{t+1}:=Y_t-\mb{e}_i+\mb{e}_j.
\]
Further, if $\chi={\sf tail}$,
then we keep both, i.e., $X_{t+1}:=x$ and $Y_{t+1}:=y$.

Note that chains $\{X_t\}_{t\in \Z_+}$ and $\{Y_t\}_{t \in \Z_+}$ are the lazy versions of the original chain.
Let us define the discrepancy
\[
D(x, y):=\frac{1}{2}\sum_{i=1}^Q|x_i-y_i|.
\]
We claim that if $X_t\neq Y_t$,
then $D(X_t, Y_t)$ increases or decreases nearly equally likely:
First, in {\bf (1a)}, the discrepancy changes by $1$ or $-1$ with equal probability
\[
\frac{1}{Q^2}\#\{x>y>0\}\cdot\#\{y>x>0\}.
\]
and, in {\bf (1b)},
the discrepancy does not change almost surely.
Next, in {\bf (2)},
the discrepancy does not increase almost surely.
Finally, in {\bf (3)},
the discrepancy does not change almost surely.

Let us assume that $x \neq y$.
Note that the following holds:
\[
\left\{x>y>0\right\}=\emptyset\quad \text{and}\quad \left\{y>x>0\right\}=\emptyset \implies
\left\{x>y=0\right\}\neq \emptyset \quad \text{and} \quad \left\{y>x=0\right\}\neq \emptyset,
\]
further that
\[
\left\{x>y>0\right\}\neq \emptyset \quad \text{and} \quad \left\{y>x>0\right\}=\emptyset \implies \left\{y>x=0\right\}\neq \emptyset.
\]
These follow since $\sum_{i=1}^Qx_i=n$ by definition of $B$.
In these two cases,
the {\bf (2)} occurs with positive probability, and the discrepancy decreases by $1$ or does not decrease with probability at least $1/Q^2$ respectively.

Summarizing all the above cases,
as far as $X_t\neq Y_t$,
the discrepancy $D_t:=D(X_t, Y_t)$ satisfies that almost surely
\[
\Eb[D_{t+1}-D_t\mid X_t, Y_t]\le 0 \quad \text{and} \quad 
\Var(D_{t+1}-D_t\mid X_t, Y_t)\ge \frac{1}{4Q^2}.
\]
Let $\tau:=\inf\{t \ge 0 \ : X_t=Y_t\}$.
There exists a constant $C_Q$ such that $\Eb_{x, y}\tau \le C_Q n^2$ for all $x, y \in B$ \cite[Proposition 17.20]{LP}.
Therefore the Markov inequality implies the claim.
\qed

\begin{remark}\label{Eq:nocutoff}
For each fixed $Q\ge 2$,
the order of the bound we obtain is optimal:
If we consider $Q=2$ and $n \ge 1$,
the chain is a symmetric nearest neighbor random walk on the path on $\{0, 1, \dots, n\}$ reflected at the extremes.
This chain does not exhibit cutoff (cf.\ \cite[Chapter 18]{LP}). 
\end{remark}

\subsection*{Acknowledgments}
The author thanks Professor Kaito Fujii for his suggestion on the chain in Appendix \ref{Sec:appendix} and for many helpful discussions and references, and Professor Makiko Sasada for her introduction to the generalized exclusion processes.
This work is partially supported by 
JSPS Grant-in-Aid for Scientific Research JP24K06711.

\bibliographystyle{alpha}
\bibliography{congestion}

\end{document}